\RequirePackage{etoolbox}
\csdef{input@path}{{style/}{graphics/}}
\documentclass[aos,MSNbibl,seceqn,dvips]{arximspdf}

\doi{10.1214/14-AOS1284}% Updated by VTEXPTS2LaTeX.exe, 17.12.2014
\volume{43}
\issue{2}
\pubyear{2015}
\firstpage{467}
\lastpage{500}
\docsubty{FLA}

\makeatletter
\renewcommand{\sharp}{\#}
\newcommand{\tr}{\operatorname{tr}}
\newcommand{\MLR}{\mathrm{MLR}}
\newcommand{\rrvert}{\vert}
\newcommand{\llvert}{\vert}
\newtheorem{teo}{Theorem}[section]
\newtheorem{prop}{Proposition}
\newproclaim{rmk}{Remark}
\newproclaim{ass}{Assumption}
\newcommand{\bbA}{{\mathbf A}}
\newcommand{\bbB}{{\mathbf B}}
\newcommand{\bbE}{{\mathbf E}}
\newcommand{\bbf}{{\mathbf f}}
\newcommand{\bbP}{{\mathbf P}}
\newcommand{\bbH}{{\mathbf H}}
\newcommand{\bbI}{{\mathbf I}}
\newcommand{\bbn}{{\mathbf n}}
\newcommand{\bbQ}{{\mathbf Q}}
\newcommand{\bbr}{{\mathbf r}}
\newcommand{\bbS}{{\mathbf S}}
\newcommand{\bbT}{{\mathbf T}}
\newcommand{\bbu}{{\mathbf u}}
\newcommand{\bbV}{{\mathbf V}}
\newcommand{\bbW}{{\mathbf W}}
\newcommand{\bbX}{{\mathbf X}}
\newcommand{\bbY}{{\mathbf Y}}
\newcommand{\bbz}{{\mathbf z}}
\newcommand{\bbL}{{\mathbf L}}
\newcommand{\bbx}{{\mathbf x}}
\newcommand{\bby}{{\mathbf y}}
\makeatother

\begin{document}
\begin{frontmatter}

\title{Independence test for high dimensional data based on
regularized canonical correlation coefficients}
\runtitle{Regularized canonical correlation coefficients}

\begin{aug}
% Corresponding author: Guangming Pan - gmpan@ntu.edu.sg% Updated by
%VTEXPTS2LaTeX.exe, 17.12.2014 10:42
\author[A]{\fnms{Yanrong}~\snm{Yang}\ead[label=e1]{yanrong.yang@monash.edu}}
\and
\author[B]{\fnms{Guangming}~\snm{Pan}\corref{}\thanksref{T1}\ead[label=e2]{GMPAN@ntu.edu.sg}} %\ead[label=u1,url]{http://www.foo.com}}
\runauthor{Y. Yang and G. Pan}
\affiliation{Monash University and Nanyang Technological University}
\address[A]{Department of Econometrics\\
\quad and Business Statistics\\
Monash University\\
Caulfield East, Victoria 3145\\
Australia\\
\printead{e1}}
\address[B]{School of Physical\\
\quad and Mathematical Sciences\\
Nanyang Technological University\\
Singapore\\
\printead{e2}}
\end{aug}
\thankstext{T1}{Supported by Ministry of Education in Singapore under
Tier 2 Grant ARC 14/11 and Tier 1 Grant RG25/14.}

% HISTORY:
%
\received{\smonth{5} \syear{2012}}% Updated by VTEXPTS2LaTeX.exe,
%17.12.2014 10:42
%
\revised{\smonth{6} \syear{2014}}% Updated by VTEXPTS2LaTeX.exe,
%17.12.2014 10:42

% ABSTRACT
%
\begin{abstract}
This paper proposes a new statistic to test independence between two
high dimensional
random vectors ${\mathbf X}\dvtx p_1\times1$ and ${\mathbf Y}\dvtx  p_2\times1$. The
proposed statistic is
based on the sum of regularized sample canonical correlation
coefficients of ${\mathbf X}$ and ${\mathbf Y}$.
The asymptotic distribution of the statistic under the null hypothesis
is established as a
corollary of general central limit theorems (CLT) for the linear
statistics of classical
and regularized sample canonical correlation coefficients when $p_1$
and $p_2$ are both
comparable to the sample size $n$.
As applications of the developed independence test, various types of
dependent structures,
such as factor models, ARCH models and a general uncorrelated but
dependent case, etc., are
investigated by simulations. As an empirical application,
cross-sectional dependence of
daily stock returns of companies between different sections in the New
York Stock Exchange (NYSE) is detected by the proposed test.
\end{abstract}

% KEYWORDS
\end{frontmatter}

%s1 #&#
\section{Introduction}\label{sec1}

A prominent feature of data collection nowadays is that the number of
variables is comparable with the sample size. This type of data poses
great challenges because traditional multivariate approaches do not
necessarily work, which were established for the case of the sample
size $n$ tending to infinity and the dimension $p$ remaining fixed (see
\cite{Ader1984}). There have been a substantial body of research work
dealing with high dimensional data, for example, \cite{BS1996,Fan2012,Huang2008,Fan2008,Baing2002,BD005}, etc.

%However, large dimensional random matrix theory always provides useful
%tools for high dimensional statistical applications.

The importance of the independence assumption for inference arises in
many aspects of multivariate analysis. For example, it is often the
case in multivariate analysis that a number of
variables can be rationally classified into several mutually exclusive
categories. When variables can be grouped in such a way, a natural
question is whether there is any significant relationship between the
groups of variables. In other words, can we claim that the groups are
mutually independent so that further statistical analysis such as
classification and testing hypothesis of equality of mean vectors and
covariance matrices could be conducted? When the dimension $p$ is
fixed, \cite{wi} used the likelihood ratio statistic to test
independence for $k$ sets of normal distributed random variables and
one may also refer to Chapter~12 of \cite{Ader1984} regarding to this point.
Relying on the asymptotic theory of sample canonical correlation
coefficients, this paper proposes a new statistic to test independence
between two high dimensional random vectors.% which are derived in
%terms of large dimensional random matrix theory.

Specifically, the aim is to test the hypothesis
%
%e1.1 #&#
\begin{eqnarray}\label{a8}
&& \mathbb{H}_0\dvtx  {\mathbf x}\mbox{ and }{\mathbf y}\mbox{ are independent};\quad\mbox{against}
\nonumber\\[-10pt]\\[-10pt]\nonumber
&& \mathbb{H}_1\dvtx  {\mathbf x}\mbox{ and }{\mathbf y}\mbox{ are dependent},
\end{eqnarray}
where ${\mathbf x}=(x_1,\ldots,x_{p_1})^T$ and ${\mathbf y}=(y_1,\ldots,y_{p_2})^T$. Without loss of generality, suppose that $p_1\leq p_2$.

It is well known that canonical correlation analysis (CCA) deals with
the correlation structure between two random vectors (see Chapter~12 of
\cite{Ader1984}). Draw $n$ independent and identically distributed
(i.i.d.) observations from these two random vectors ${\mathbf x}$ and ${\mathbf
y}$, respectively, and group them into $p_1\times n$ random matrix
$\bbX=({\mathbf x}_1,\ldots,{\mathbf x}_n)=(X_{ij})_{p_1\times n}$ and
$p_2\times n$ random\vspace*{1pt} matrix $\bbY=({\mathbf y}_1,\ldots,{\mathbf
y}_n)=(Y_{ij})_{p_2\times n}$, respectively. CCA seeks the linear
combinations ${\mathbf a}^{T}{\mathbf x}$ and ${\mathbf b}^{T}{\mathbf y}$ that are
most highly correlated, that is, to maximize
%
%e1.2 #&#
\begin{equation}
\gamma=\operatorname{Corr}\bigl({\mathbf a}^{T}{\mathbf x}, {\mathbf b}^{T}{\mathbf y}
\bigr)=\frac{{\mathbf
a}^{T}\bolds\Sigma_{{\mathbf x}{\mathbf y}}{\mathbf b}}{\sqrt{{\mathbf
a}^{T}\bolds\Sigma_{{\mathbf x}{\mathbf x}}{\mathbf a}}\sqrt{{\mathbf b}^{T}
\bolds\Sigma_{{\mathbf y}{\mathbf y}}{\mathbf b}}},
\end{equation}
where $\bolds\Sigma_{{\mathbf x}{\mathbf x}}$ and $\bolds\Sigma
_{{\mathbf y}{\mathbf y}}$ are the population covariance matrices for ${\mathbf x}$
and ${\mathbf y}$, respectively, and
$\bolds\Sigma_{{\mathbf x}{\mathbf y}}$ is the population covariance
matrix between ${\mathbf x}$ and ${\mathbf y}$. After finding the maximal
correlation $r_1$ and associated vectors ${\mathbf a}_1$ and ${\mathbf b}_1$,
CCA continues to seek a second linear combination ${\mathbf a}_2^{T}{\mathbf
x}$ and ${\mathbf b}_2^{T}{\mathbf y}$ that has the maximal correlation among
all linear combinations uncorrelated with ${\mathbf a}_1^{T}{\mathbf x}$ and
${\mathbf b}_1^{T}{\mathbf y}$. This procedure can be iterated and successive
canonical correlation coefficients $\gamma_1,\ldots,\gamma_{p_1}$
can be found.

It turns out that the population canonical correlation coefficients
$\gamma_1,\ldots,\gamma_{p_1}$ can be recast as the roots of the
determinant equation
%
%e1.3 #&#
\begin{equation}
\label{root} \det\bigl(\bolds\Sigma_{{\mathbf x}{\mathbf y}}\bolds
\Sigma_{{\mathbf
y}{\mathbf y}}^{-1}\bolds\Sigma_{{\mathbf x}{\mathbf y}}^{T}-
\gamma^{2}\bolds\Sigma_{{\mathbf x}{\mathbf x}}\bigr)=0.
\end{equation}
Regarding this point, one may refer to page $284$ of \cite{MKB}. The
roots of the determinant equation above go under many names,
because they figure equally in discriminant analysis, canonical
correlation analysis and invariant tests of linear hypotheses in
the multivariate analysis of variance. %These are standard techniques
%in multivariate statistical analysis.

Traditionally, sample covariance matrices $\hat{\bolds\Sigma
}_{{\mathbf x}{\mathbf x}}$, $\hat{\bolds\Sigma}_{{\mathbf x}{\mathbf y}}$ and
$\hat{\bolds\Sigma}_{{\mathbf y}{\mathbf y}}$ are used to replace the
corresponding population covariance matrices to solve the nonnegative
roots $\rho_1, \rho_2, \ldots, \rho_{p_1}$ to the determinant equation
\[
\det\bigl(\hat{\bolds\Sigma}_{{\mathbf x}{\mathbf y}} \hat{\bolds
\Sigma}_{{\mathbf y}{\mathbf y}}^{-1}\hat{\bolds\Sigma}_{{\mathbf x}{\mathbf
y}}^{T}-
\rho^{2}\hat{\bolds\Sigma}_{{\mathbf x}{\mathbf x}}\bigr)=0,
\]
where
\begin{eqnarray*}
\hat{\bolds\Sigma}_{{\mathbf x}{\mathbf x}}&=&\frac{1}{n}\sum
^{n}_{i=1}({\mathbf x}_i-\bar{\mathbf x}) ({\mathbf
x}_i-\bar{\mathbf x})^{T},\qquad \hat{\bolds
\Sigma}_{{\mathbf x}{\mathbf y}}=\frac{1}{n}\sum^{n}_{i=1}({
\mathbf x}_i-\bar{\mathbf x}) ({\mathbf y}_i-\bar{\mathbf
y})^{T}, %
\\
\hat{\bolds\Sigma}_{{\mathbf y}{\mathbf y}}&=&\frac{1}{n}\sum
^{n}_{i=1}({\mathbf y}_i-\bar{\mathbf y}) ({\mathbf
y}_i-\bar{\mathbf y})^{T},\qquad \bar{\mathbf x}=\frac{1}{n}\sum
^{n}_{i=1}{\mathbf x}_i,\qquad \bar{\mathbf
y}=\frac
{1}{n}\sum^{n}_{i=1}{\mathbf
y}_i.
\end{eqnarray*}
However, it is inappropriate to use these types of sample covariance
matrices to replace population covariance matrices to test (\ref{a8})
in some cases. We demonstrate such an example in Section~\ref{factor2}.

%consider two random vectors
%where $\tau$ is a random variable and $\bbe=(1,\cdots,1)$. Evidently, $
%based on the above sample covariance matrices will conclude that $
%because the above sample covariance matrices treat $\bbu$ and $\bbv$
%as $\bbx$ and $y$.

Therefore, in this paper we instead consider the nonnegative roots
$r_1,r_2,\ldots,\break r_{p_1}$ of an alternative determinant equation as follows:
%
%e1.4 #&#
\begin{equation}
\label{root1} \det\bigl(\bbA_{{\mathbf x}{\mathbf y}}\bbA_{{\mathbf y}{\mathbf y}}^{-1}
\bbA_{{\mathbf
x}{\mathbf y}}^{T}-r^{2}\bbA_{{\mathbf x}{\mathbf x}}\bigr)=0,
\end{equation}
where
\[
\bbA_{{\mathbf x}{\mathbf x}}=\frac{1}{n}\bbX\bbX^T,\qquad
\bbA_{{\mathbf y}{\mathbf y}}=\frac{1}{n}\bbY\bbY^T,\qquad
\bbA_{{\mathbf x}{\mathbf y}}=\frac{1}{n}\bbX\bbY^T.
\]
%
%because they are consistent estimators of respective population
%covariance matrices $\bolds\Sigma_{\bbx\bbx}$, $\bolds
%In this case, substituting sample covariance matrices for
%corresponding population covariance matrices in (\ref{root}), are
%consistent estimators of canonical correlation coefficients $\gamma_1,
We also call $\bbA_{{\mathbf x}{\mathbf x}}$, $\bbA_{{\mathbf y}{\mathbf y}}$ and
$\bbA_{{\mathbf x}{\mathbf y}}$ sample covariance matrices, as in the random
matrix community. However, whichever sample covariance matrices are
used they are not consistent estimators of population covariance
matrices, which is called the ``curse of dimensionality,'' when the
dimensions $p_1$ and $p_2$ are both comparable to the sample size $n$.
As a consequence, it is conceivable that the classical likelihood ratio
statistic (see \cite{wi} and \cite{Ader1984}) does not work well in
the high dimensional case (in fact, it is not well defined and we will
discuss this point in the later section).

Moreover, from (\ref{root1}), when $p_1<n,p_2<n$, one can see that
$r_1^2, r_2^2, \ldots, r_{p_1}^2$ are the eigenvalues of the matrix
%
%e1.5 #&#
\begin{equation}
\label{a3*} \bbS_{{\mathbf x}{\mathbf y}}=\bbA_{{\mathbf x}{\mathbf x}}^{-1}
\bbA_{{\mathbf x}{\mathbf
y}}\bbA_{{\mathbf y}{\mathbf y}}^{-1}\bbA_{{\mathbf x}{\mathbf y}}^{T}.
\end{equation}
Evidently, $\bbA_{{\mathbf x}{\mathbf x}}^{-1}$ and $\bbA_{{\mathbf y}{\mathbf
y}}^{-1}$ do not exist when $p_1>n$ and $p_2>n$. For this reason, we
also consider the eigenvalues of the regularized matrix
%
%e1.6 #&#
\begin{equation}
\label{a3} \bbT_{{\mathbf x}{\mathbf y}}=\bbA_{t{\mathbf x}}^{-1}
\bbA_{{\mathbf x}{\mathbf y}}\bbA_{{\mathbf y}{\mathbf y}}^{-}\bbA_{{\mathbf x}{\mathbf y}}^{T},
\end{equation}
where $\bbA_{t{\mathbf x}}^{-1}=(\frac{1}{n}\bbX\bbX^{T}+t\bbI
_{p_1})^{-1}$, $t$ is a positive constant number and $\bbI_{p_1}$ is a
$p_1\times p_1$ identity matrix, and $\bbA_{{\mathbf y}{\mathbf y}}^{-}$
denotes\vspace*{1pt} the Moore--Penrose pseudoinverse matrix of~$\bbA_{{\mathbf y}{\mathbf
y}}$. Hence, $\bbT_{{\mathbf x}{\mathbf y}}$ is well defined even in the case
of $p_1,p_2\geq n$.
Moreover, $\bbT_{{\mathbf x}{\mathbf y}}$ reduces to $\bbS_{{\mathbf x}{\mathbf y}}$
when $p_1,p_2$ are both smaller than $n$ and $t=0$.
%The linear spectral statistic of the matrix $\bbT_{\bbx\bby}$ can deal
%with high dimensional random vectors $\bbx\dvtx  p_1\times1$ and $\bby\dvtx  p_2
%matrix $(\frac{1}{n}\bbY\bbY^{T})^{-}$ is equivalent to the classical
%inverse matrix $(\frac{1}{n}\bbY\bbY^{T})^{-1}$ and take $t=0$ in $

We now look at CCA from another perspective. The original random
vectors ${\mathbf x}$ and ${\mathbf y}$ can be transformed into new random
vectors $\bolds{\xi}$ and $\bolds{\eta}$ as
%
%e1.7 #&#
\begin{eqnarray}
\label{02} \pmatrix{ {\mathbf x}
\cr
{\mathbf y}}\rightarrow\pmatrix{ \bolds{\xi}
\cr
\bolds{\eta}}= \pmatrix{ \mathcal{A}' & \bolds{0}
\cr
\bolds{0} & \mathcal{B}'} \pmatrix{ {\mathbf x}
\cr
{\mathbf y}}
\end{eqnarray}
such that
%
%e1.8 #&#
\begin{eqnarray}
\label{01} \pmatrix{ \mathcal{A}' & \bolds{0}
\cr
\bolds{0} &
\mathcal{B}'}\pmatrix{ \bolds{\Sigma}_{{\mathbf x}{\mathbf x}} &
\bolds{\Sigma}_{{\mathbf
x}{\mathbf y}}
\cr
\bolds{\Sigma}_{{\mathbf y}{\mathbf x}} &
\bolds{\Sigma}_{{\mathbf
y}{\mathbf y}}} \pmatrix{ \mathcal{A} & \bolds{0}
\cr
\bolds{0} & \mathcal{B}}= \pmatrix{ \bbI_{p_1} & \mathcal{P}
\cr
\mathcal{P}' & \bbI_{p_2}},
\end{eqnarray}
where $\mathcal{P}=(\mathcal{P}_1,\bolds{0})$, $\mathcal
{P}_1=\operatorname{diag}(\gamma_1, \ldots, \gamma_{p_1})$ and
$
\mathcal{A}=\bolds{\Sigma}_{{\mathbf x}{\mathbf x}}^{-1/2}\bbQ_1, \mathcal
{B}=\bolds{\Sigma}_{{\mathbf y}{\mathbf y}}^{-1/2}\bbQ_2$,
with $\bbQ_1\dvtx  p_1\times p_1$ and $\bbQ_2\dvtx  p_2\times p_2$ being
orthogonal matrices satisfying
\begin{eqnarray*}
&&\bolds{\Sigma}_{{\mathbf x}{\mathbf x}}^{-1/2}\bolds{\Sigma
}_{{\mathbf x}{\mathbf y}}\bolds{\Sigma}_{{\mathbf y}{\mathbf y}}^{-1/2}=\bbQ
_1\mathcal{P}\bbQ_2.
\end{eqnarray*}
Hence, testing independence between ${\mathbf x}$ and ${\mathbf y}$ is
equivalent to testing independence between $\bolds{\xi}$ and
$\bolds{\eta}$. The covariance between $\bolds{\xi}$ and
$\bolds{\eta}$ has the following simple expression
%
%e1.9 #&#
\begin{eqnarray}
\operatorname{Var}\pmatrix{ \bolds{\xi}
\cr
\bolds{\eta}}= \pmatrix{
\bbI_{p_1} & \mathcal{P}
\cr
\mathcal{P}' &
\bbI_{p_2}}.
\end{eqnarray}
In view of this, if the joint distribution of ${\mathbf x}$ and ${\mathbf y}$
is Gaussian, independence between ${\mathbf x}$ and ${\mathbf y}$ is equivalent
to asserting that the
population canonical correlations all vanish: $\gamma_1=\cdots=
\gamma_{p_1}=0$. Details can be referred to Chapter~11 of \cite
{Fujikoshi2010}.
A~natural criteria for this test should be $\sum^{p_1}_{i=1}\gamma_i^2$.

As pointed out, $r_i$ is not a consistent estimator of the
corresponding population
version $\gamma_i$ in the high dimensional case. However,
fortunately, the classical sample canonical correlation coefficients
$r_1,r_2,\ldots,r_{p_1}$ or its regularized analogous still contain important
information so that hypothesis testing for (\ref{a8}) is possible
although the classical likelihood ratio statistic does not work well
in the high dimensional case. This is due to the fact that the limits
of the empirical spectral distribution (ESD) of $r_1,\ldots,r_{p_1}$
under the null and the alternative could be different so that we may
use it to distinguish dependence from independence (one may see the
next section).
Our approach essentially makes use of the integral of functions with
respect to the ESD of canonical correlation coefficients.
The proposed statistic turns out a trace of the corresponding matrices,
that is, $\sum^{p_1}_{i=1}r_i^2$. In order to apply it to conduct
tests, we further propose
two modified statistics by either dividing the total samples into two
groups or estimating the population covariance matrix of ${\mathbf x}$ in a
framework of sparsity.

In addition to proposing a statistic for testing (\ref{a8}), another
contribution of this paper is to establish the limit of the ESD of
regularized sample canonical correlation coefficients and central limit
theorems (CLT) of linear functionals of the classical and regularized
sample canonical correlation coefficients $r_1,r_2,\ldots,r_{p_1}$,
respectively. %i.e. the CLT of $\sum^{p_1}_{i=1}\phi(r_i^2)$, where $
%on some specific subsets.
This is of an independent interest in its own right in addition to
providing asymptotic distributions for the proposed statistics.

To derive the CLT for linear spectral statistics of classical and
regularized sample canonical correlation coefficients, the strategy is
to first establish the CLT under the Gaussian case, that is, the
entries of $\bbX$ are Gaussian distributed. In the Gaussian case, the
CLT for linear spectral statistics of the matrix $\bbS_{{\mathbf x}{\mathbf
y}}$ can be linked to that of an $F$-matrix, which has been
investigated in \cite{zsr}. We then extend the CLT to general
distributions by bounding the difference between the characteristic\vadjust{\goodbreak}
functions of the respective linear spectral statistics of $\bbS_{{\mathbf
x}{\mathbf y}}$ under the Gaussian case and non-Gaussian case. To bound
such a difference and handle the inverse of a random matrix, we use an
interpolation approach and a smooth cutoff function. The approach of
developing the CLT for linear spectral statistics of the matrix $\bbT
_{{\mathbf x}{\mathbf y}}$ is similar to that for $\bbS_{{\mathbf x}{\mathbf y},}$
except we first have to develop CLT of perturbed sample covariance
matrices in the supplement material \cite{supp} for establishing CLT of the matrix
$\bbT_{{\mathbf x}{\mathbf y}}$ when the entries of $\bbX$ are Gaussian.

%Due to existence of some unknown parameter involving $\bolds{
%matrix $\bbT_{\bbx\bby}$, we provide two methods to deal with it. One
%is to estimate it for a class of sparse matrices and the other one is
%to eliminate it by a dividing-sample method.

Here, we would point out some works on canonical correlation
coefficients under the high dimensional scenario. In the high
dimensional case, \cite{WK}
investigated the limit of the ESD of the classical sample canonical
correlation coefficients $r_1,r_2,\ldots,r_{p_1}$ and \cite{Jo}
established the Tracy--Widom law of the maximum of sample correlation
coefficients when $\bbA_{{\mathbf x}{\mathbf x}}$ and $\bbA_{{\mathbf y}{\mathbf y}}$
are Wishart matrices and ${\mathbf x}$, ${\mathbf y}$ are independent.

%Under the assumptions of Gaussian entries and mutual independence of
%the matrices $\bbX$ and $\bbY$, each eigenvalue of $\bbS_{\bbx\bby}$
%is equivalent to a function of the corresponding eigenvalue of an
%$F$-matrix. \cite{zsr} has provided the CLT for an $F$-matrix. So the
%CLT for canonical correlation coefficients in the Gaussian case can be
%derived. Moreover, an `interpolation' trick that is proposed in
%the same as that in the Gaussian case.

%As is known to all, the uncorrelated relationship is equivalent to the
%independent relationship only when the tested random vectors are
%Gaussian distributed. For nonGaussian case, we can utilize the
%correlation between the high powers of the elements of random vectors
%to test the dependence of the initial random vectors. This is also the
%efficiency of developing the CLT in nonGaussian case for independence
%test.

%As an illustration of statistical applications of the CLT for sample
%canonical correlation coefficients, we establish the application of
%the proposed independence test for no linear relationship test and
%cointegration test in high dimensional case.

The remainder of the paper is organized as follows. Section~\ref{sec2}
proposes a new test statistic for (\ref{a8}) based on large
dimensional random matrix theory and contains the main results.
%%Section $3$ proposes based on the CLT given in Section $2$.
Two modified statistics are further provided in Section~\ref{sec3} and
estimators for some unknown parameters in the asymptotic mean and
variance for the asymptotic distribution are also proposed.
Section~\ref{sec4} provides the powers of the test statistics. Two examples as
statistical inference of independence test are explored in Section~\ref{sec5}.
Simulation results for several kinds of dependent structures are
provided in Section~\ref{sec6}. An empirical analysis of cross-sectional
dependence of daily stock returns of companies from two different
sections in New York Stock Exchange (NYSE) is investigated by the
proposed independence test in Section~\ref{sec7}. Some useful lemmas and
proofs of all theorems and Propositions~\ref{pr4}--\ref{pr5} are relegated to Appendix A while
one theorem about the CLT of a sample covariance matrix plus a
perturbation matrix is provided in Appendix B. All appendices are given
in the supplementary material \cite{supp}.

%s2 #&#
\section{Methodology and theory}\label{sec2}

Throughout this paper, we make the following assumptions.

%as1 #&#
\begin{ass}\label{ass1}
$p_1=p_1(n)$ and $ p_2=p_2(n)$ with $\frac
{p_1}{n}\rightarrow c_1$ and $ \frac{p_2}{n}\rightarrow c_2$, $c_1,
c_2\in(0,1)$, as $n\rightarrow\infty$.
\end{ass}

%as2 #&#
\begin{ass}\label{ass2}
$p_1=p_1(n)$ and $ p_2=p_2(n)$ with $\frac
{p_1}{n}\rightarrow c'_1$ and $ \frac{p_2}{n}\rightarrow c'_2$,
$c'_1\in(0,+\infty)$ and $c_2'\in(0,+\infty)$, as
$n\rightarrow\infty$.
\end{ass}

%as3 #&#
\begin{ass}\label{ass3}
$\bbX=(X_{ij})_{i,j=1}^{p_1,n}$ and $\bbY
=(Y_{ij})_{i,j=1}^{p_2,n}$ satisfy $\bbX=\bolds{\Sigma}_{{\mathbf
x}{\mathbf x}}^{1/2}\bbW$ and $\bbY=\bolds{\Sigma}_{{\mathbf y}{\mathbf
y}}^{1/2}\bbV$, where $\bbW=({\mathbf w}_1,\ldots,{\mathbf
w}_n)=(W_{ij})_{i,j=1}^{p_1,n}$ consists of i.i.d. real random variables
$\{W_{ij}\}$ with $EW_{11}=0$ and
$E\llvert W_{11}\rrvert ^2=1$; $\bbV=({\mathbf v}_1,\ldots,{\mathbf
v}_n)=(V_{ij})_{i,j=1}^{p_2,n}$ consists of i.i.d. real random
variables with $EV_{11}=0$ and
$E\llvert V_{11}\rrvert ^2=1$; $\bolds{\Sigma}_{{\mathbf x}{\mathbf
x}}^{1/2}$ and
$\bolds{\Sigma}_{{\mathbf y}{\mathbf y}}^{1/2}$ are Hermitian square
roots of positive definite matrices $\bolds{\Sigma}_{{\mathbf x}{\mathbf
x}}$ and $\bolds{\Sigma}_{{\mathbf y}{\mathbf y}}$, respectively, so
that $(\bolds{\Sigma}_{{\mathbf x}{\mathbf x}}^{1/2})^2=\bolds
{\Sigma}_{{\mathbf x}{\mathbf x}}$ and $(\bolds{\Sigma}_{{\mathbf y}{\mathbf
y}}^{1/2})^2=\bolds{\Sigma}_{{\mathbf y}{\mathbf y}}$.
\end{ass}

%as4 #&#
\begin{ass}\label{ass4}
$F^{\bolds{\Sigma}_{{\mathbf x}{\mathbf
x}}}\stackrel{D}\rightarrow H$, a proper cumulative distribution function.
\end{ass}

%re1 #&#
\begin{rmk}
By the definition of the matrix $\bbS_{{\mathbf x}{\mathbf y}}$, the classical
canonical correlation coefficients between ${\mathbf x}$ and ${\mathbf y}$ are
the same as those between ${\mathbf w}$ and ${\mathbf v}$ when ${\mathbf w}$ and $\{
{\mathbf w}_i\}$ are i.i.d., and ${\mathbf v}$ and $\{{\mathbf v}_i\}$ are i.i.d.
\end{rmk}

We now introduce some results from random matrix theory. Denote the ESD
of any $n\times n$ matrix $\bbA$ with real eigenvalues
$\mu_1\leq\mu_2\leq\cdots\leq\mu_n$ by
%
%e2.1 #&#
\begin{equation}
\label{a12} F^{\bbA}(x)=\frac{1}{n}\#\{i\dvtx  \mu_i\leq
x\},
\end{equation}
where $\#\{\cdots\}$ denotes the cardinality of the set $\{\cdots\}$.
%%From the determinant equation (\ref{root1}), it can be seen that $
%the matrix $\bbS_{\bbx\bby}$. Hence analyzing the ESD, $F^{\bbS_{\bbx
%the ESD of $r_1, r_2, \cdots, r_{p_1}$. Here the ESD of $r_1, r_2,
%$r_i$. %Thus the limit of the ESD, $F^{\bbS_{\bbx\bby}}(x)$, of $\bbS_{
%canonical correlation coefficients $r_1, r_2, \cdots, r_{p_1}$, the
%roots to (\ref{root1}).
%F_n(x)=\frac{1}{p_1}\#\{i\dvtx  r_i\leq x\}.

When the two random vectors ${\mathbf x}$ and ${\mathbf y}$ are independent and
each of them consists of i.i.d. Gaussian random
variables, under Assumptions \ref{ass1}~and~\ref{ass3}, \cite{WK} proved that the
empirical measure of the classical sample canonical correlation
coefficients $r_1, r_2, \ldots, r_{p_1}$ converges in
probability to a fixed distribution whose density is given by
%
%e2.2 #&#
\begin{equation}
\label{lsd} \rho(x)=\frac{\sqrt{(x-L_1)(x+L_1)(L_2-x)(L_2+x)}}{\pi
c_1x(1-x)(1+x)},\qquad x\in[L_1,L_2],
\end{equation}
and atoms size of $\max(0,(1-c_2)/c_1)$ at zero and size $\max
(0,1-(1-c_2)/c_1)$ at unity where $L_1=\llvert \sqrt{c_2-c_2c_1}-\sqrt
{c_1-c_1c_2}\rrvert $ and $L_2=\llvert \sqrt{c_2-c_2c_1}+\break  \sqrt
{c_1-c_1c_2}\rrvert $. Here,
the empirical measure of $r_1, r_2, \ldots, r_{p_1}$ is defined as in
(\ref{a12}) with $\mu_i$ replaced by $r_i$.

\cite{yp} proved that (\ref{lsd}) also holds for classical
sample canonical correlation coefficients when the entries of ${\mathbf x}$
and ${\mathbf y}$ are not necessarily Gaussian distributed. For easy
reference, we state the result in the following proposition.

%pr1 #&#
\begin{prop}\label{theo1}
In addition to Assumptions \ref{ass1} and \ref{ass3}, suppose that\break 
%$(d)$ $\bbS_{\bbx\bby}=\bbA_{\bbx\bbx}^{-1}\bbA_{\bbx\bby}\bbA_{\bby
$\{X_{ij}, 1\leq i\leq p_1, 1\leq j\leq n\}$ and $\{Y_{ij}, 1\leq i\leq
p_2, 1\leq j\leq n\}$ are independent. Then the empirical measure of
$r_1,r_2,\ldots,r_{p_1}$ converges almost surely to a fixed
distribution function whose density is given by (\ref{lsd}).
\end{prop}

Under Assumptions \ref{ass2}--\ref{ass4}, instead of $F^{\bbS_{{\mathbf x}{\mathbf y}}}$, we
analyze the ESD, $F^{\bbT_{{\mathbf x}{\mathbf y}}}$, of the regularized
random matrix $\bbT_{{\mathbf x}{\mathbf y}}$ given in (\ref{a3}). To this
end, define the Stieltjes
transform of any distribution function $G(x)$ by
\[
m_G =\int\frac{1}{x-z}\,dG(x),\qquad z\in\mathbb{C}^+=\{z\in\mathbb
{C},\Im z>0\}, %
\]
where $\Im z$ denotes the imaginary part of the complex number $z$.

It turns out that the limit of the empirical spectral distribution
(LSD) of $\bbT_{{\mathbf x}{\mathbf y}}$ is connected to the LSD of $\bbS
_1\bbS^{-1}_{2t}$ defined below. Let
\begin{eqnarray*}
&&\bbS_1=\frac{1}{p_2}\sum^{p_2}_{k=1}{
\mathbf w}_k{\mathbf w}_k^{T},\qquad
\bbS_{2t}=\frac{1}{n-p_2}\sum^{n}_{k=p_2+1}{
\mathbf w}_k{\mathbf w}_k^{T}+t\frac{n}{n-p_2}
\bolds{\Sigma}_{{\mathbf x}{\mathbf x}}^{-1},
\\
&&y_1=\frac{c_1'}{c_2'}, \qquad y_2=\frac{c_1'}{1-c_2'}.
\end{eqnarray*}
In the definition of $\bbS_{2t}$, we require $n>p_2$. The LSD of $\bbS
_{2t}$ and its Stieltjes transform are denoted by $F_{y_2t}$ and
$m_{y_2t}(z)$, respectively. Under Assumptions~\mbox{\ref{ass2}--\ref{ass4}}, from~\cite{SB1995}
and \cite{Pan}, $m_{y_2t}(z)$ is the unique solution in $\mathbb
{C}^+$ to
%
%e2.3 #&#
\begin{equation}
\label{b6} m_{y_2t}(z)=m_{H_t} \biggl(z-\frac{1}{1+y_2m_{y_2t}(z)}
\biggr),
\end{equation}
where\vspace*{1pt} $m_{H_t}(z)$ denotes the Stieltjes transform of the LSD of the
matrix $t\frac{n}{n-p_2}\bolds{\Sigma}_{{\mathbf x}{\mathbf x}}^{-1}$
(one may also see (1.4) in the supplement material \cite{supp}). Let $\bbn
=(n_1,n_2)$ and ${\mathbf y}=(y_1,y_2)$ with $n_1=p_1$ and $n_2=n-p_2$. The
Stieltjes transforms of the ESD and LSD of the matrix $\bbS_1\bbS
_{2t}^{-1}$ are denoted by $m_{\bbn}(z)$ and $m_{{\mathbf y}}(z)$,
respectively, while those of the ESD and LSD of the matrix $\bbW
_1^{T}\bbS_{2t}^{-1}\bbW_1$ are denoted by $\underline{m}_{\bbn
}(z)$ and $\underline{m}_{{\mathbf y}}(z)$, respectively, where $\bbW
_1=({\mathbf w}_1, {\mathbf w}_2, \ldots, {\mathbf w}_{p_2})$. Observe that the
spectral of $\bbS_1\bbS^{-1}_{2t}$ and $\bbW_1^{T}\bbS
_{2t}^{-1}\bbW_1$ are the same except zero eigenvalues and this leads to
%
%e2.4 #&#
\begin{equation}
\label{w5} \underline{m}_{{\mathbf y}}(z)=-\frac{1-y_1}{z}+y_1m_{{\mathbf y}}(z).
\end{equation}

We are now in a position to state the LSD of $\bbT_{{\mathbf x}{\mathbf y}}$.

%th2.1 #&#
\begin{teo}\label{teo2}
In addition to Assumptions \ref{ass2}--\ref{ass4}, suppose that $\{X_{ij}, 1\leq i\leq
p_1, 1\leq j\leq n\}$ and $\{Y_{ij}, 1\leq i\leq p_2, 1\leq j\leq n\}$
are independent.

\begin{longlist}[(a)]
\item[(a)] If $c_2'\in(0,1)$, then the ESD, $F^{\bbT_{{\mathbf x}{\mathbf
y}}}(\lambda)$, converges almost surely to a fixed distribution
$\widetilde F(\frac{\lambda}{q(1-\lambda)})$ with $q=\frac
{c_2'}{1-c_2'}$ where $\widetilde F(\lambda)$ is a nonrandom
distribution and its Stieltjes transform $m_{{\mathbf y}}(z)$ is the unique
solution in $\mathbb{C}^+$ to
%
%e2.5 #&#
\begin{equation}
\label{bb2}
m_{{\mathbf y}}(z) =-\int\frac{dF_{y_2t}(1/\lambda)}{\lambda
(1-y_1-y_1zm_{{\mathbf y}}(z))-z}. %\frac{1}{z}m_{\widetilde F}(\frac{1}{z})=
%iy_{1}'z}\oint_{\left\vert\zeta\right\vert=\rho}\log\Big(1-z
%+\zeta^{-2}zy'_1m_{y_2'}(\frac{1}{\zeta})\Big)\,d\zeta,
\end{equation}
%
%for any $\rho\in(0,(1-\sqrt{1-c_2'})^2+t)$ with; where $m_{y_1'}$
%is the Stieltjes transform of the LSD of $\frac{1}{n-p_2}\bbW_2
%index $y_2'$, i.e.
%m_{y_2'}(\zeta)=\frac{1-c_2}{t}m_{H}\Big(\frac{1-c_2}{t}\big(z-
%with $m_H$ being the Stieltjes transform of the LSD for the matrix $

\item[(b)] If\vspace*{1pt} $c_2'\in[1,\infty)$, then $F^{\bbT_{{\mathbf x}{\mathbf
y}}}(\lambda)$, converges almost surely to a fixed distribution
$\widetilde{G}(\frac{t}{1-\lambda}-t)$ where $\widetilde G(\lambda)$ is a\vspace*{1pt}
nonrandom distribution and its Stieltjes transform satisfies the equation
%
%e2.6 #&#
\begin{equation}
m_{\widetilde G}(z) =\int\frac{dH(\lambda)}{\lambda
(1-c_1'-c_1'zm_{\widetilde G}(z))-z}.
\end{equation}
\end{longlist}
\end{teo}

%re2 #&#
\begin{rmk}
Indeed, taking $t=0$ in (\ref{bb2}) recovers \cite{WK}'s result (one
may refer to the result of $F$ matrix in \cite{BS2009}).
\end{rmk}
%
%the limiting spectral distribution (LSD) is

Let us now introduce the test statistic. Under Assumptions \ref{ass1}~and~\ref{ass3}, behind our test statistic is the observation that the
limit of $F^{\bbS_{{\mathbf x}{\mathbf y}}}(x)$ can be obtained from~(\ref
{lsd}) when ${\mathbf x}$ and ${\mathbf y}$ are independent, while the limit of
$F^{\bbS_{{\mathbf x}{\mathbf y}}}(x)$ could be different from~(\ref{lsd})
when ${\mathbf x}$ and ${\mathbf y}$ have correlation. For example, if ${\mathbf
y}=\Sigma_1{\mathbf w}$ and ${\mathbf x}=\Sigma_2{\mathbf w}$ with $p_1=p_2$ and
both $\Sigma_1$ and $\Sigma_2$ being invertible, then
\[
\bbS_{{\mathbf x}{\mathbf y}}=\bbI, %
\]
which implies that the limit of $F^{\bbS_{{\mathbf x}{\mathbf y}}}(x)$ is a
degenerate distribution. This suggests that we may make use of $F^{\bbS
_{{\mathbf x}{\mathbf y}}}(x)$ to construct a test statistic. Thus, we consider
the following statistic:
%
%e2.7 #&#
\begin{equation}
\label{b1} \int\phi(x) \,dF^{\bbS_{{\mathbf x}{\mathbf y}}}(x)=\frac{1}{p_1}\sum
^{p_1}_{i=1}\phi\bigl(r_i^{2}
\bigr).
\end{equation}
A perplexing problem is how to choose an appropriate function $\phi
(x)$. For simplicity, we choose $\phi(x)=x$ in this work. That is, our
statistic is
%
%e2.8 #&#
\begin{equation}
\label{b3} S_n=\int x \,dF^{\bbS_{{\mathbf x}{\mathbf y}}}(x)=\frac{1}{p_1}
\sum^{p_1}_{i=1}r_i^{2}.
\end{equation}
Indeed, extensive simulations based on Theorems \ref{teo1} and \ref
{teo3} below have been conducted to help select an appropriate function
$\phi(x)$. We find that other functions such as $\phi(x)=x^2$ does
not have an advantage over $\phi(x)=x$.

In the classical CCA, the maximum likelihood ratio test statistic for
(\ref{a8}) with fixed dimensions is
%
%e2.9 #&#
\begin{equation}
\label{b2} \MLR_n=\sum^{p_1}_{i=1}
\log\bigl(1-r_i^2\bigr)
\end{equation}
(see \cite{wi} and \cite{Ader1984}). That is, $\phi(x)$ in (\ref
{b1}) takes $\log(1-x)$. Note that the density $\rho(x)$ has atom
size of $\max(0,1-(1-c_2)/c_1)$ at unity by~(\ref{lsd}).
Thus, the normalized statistic $\MLR_n$ is not well defined when
\mbox{$c_1+c_2>1$} [because $\int\log(1-x^2)\rho(x)\,dx$ is not meaningful].
In addition, even when $c_1+c_2\leq1$, the right end point of $\rho
(x)$, $L_2$, can be equal to one so that some sample correlation
coefficients $r_i$ are close to one, for example, $L_2=1$ when
$c_1=c_2=1/2$. This in turns causes a big value of the corresponding
$\log(1-r_i^2)$. Therefore, $\MLR_n$ is not stable and this phenomenon
is also confirmed by our simulations. %Also it is easily seen that
%$H=1$ when $c_1=c_2=1/2$ or $c_1=1/4$ and $c_2=3/4$ (indeed, there are
%many other possibilities of $c_1$ and $c_2$ so that $H=1$). Therefore
%the normalized statistic $\MLR_n$ will have singular values.

Under Assumptions \ref{ass2}--\ref{ass4}, we substitute $\bbT_{{\mathbf x}{\mathbf y}}$ for $\bbS
_{{\mathbf x}{\mathbf y}}$ and use the statistic
%
%e2.10 #&#
\begin{equation}
\label{bb1} T_n=\int x\,dF^{\bbT_{{\mathbf x}{\mathbf y}}}(x).
\end{equation}

We next establish the CLTs of the statistics $(\ref{b1})$ and $(\ref
{bb1})$. To this end, write
%
%e2.11 #&#
\begin{equation}
\label{c1} G^{(1)}_{p_1,p_2}(\lambda)=p_1
\bigl(F^{\bbS_{{\mathbf x}{\mathbf y}}}(\lambda)-F^{c_{1n},c_{2n}}(\lambda)
\bigr),
\end{equation}
and
%
%e2.12 #&#
\begin{equation}
\label{c2} G^{(2)}_{p_1,p_2}(\lambda)=p_1
\bigl(F^{\bbT_{{\mathbf x}{\mathbf y}}}(\lambda
)-F^{c'_{1n},c'_{2n}}(\lambda) \bigr),
\end{equation}
where $F^{c_{1n},c_{2n}}(\lambda)$ and
$F^{c'_{1n},c'_{2n}}(\lambda)$ are obtained from
$F^{c_1,c_2}(\lambda)$ and $F^{c'_1,c'_2}(\lambda)$ with
$c_1,c_2$, $c'_1,c'_2$ and $H$ replaced by $c_{1n}=\frac
{p_1}{n},c_{2n}=\frac{p_2}{n}$, $c'_{1n}=\frac
{p_1}{n},c'_{2n}=\frac{p_2}{n}$ and $F^{\bolds{\Sigma}_{{\mathbf
x}{\mathbf x}}}$, respectively; $F^{c_1,c_2}(\lambda)$ and
$F^{c'_1,c'_2}(\lambda)$ are the limiting spectral distributions
of the matrices $\bbS_{{\mathbf x}{\mathbf y}}$ and $\bbT_{{\mathbf x}{\mathbf y}}$,
respectively. The density of $F^{c_1,c_2}(\lambda)$ can be obtained
from $\rho(x)$ in (\ref{lsd}) while the density of
$F^{c'_1,c'_2}(\lambda)$ can be recovered from (\ref{bb2}). We
renormalize $(\ref{b1})$ and $(\ref{bb1})$ as
%
%e2.13 #&#
\begin{equation}
\label{cen} \qquad\int\phi(\lambda)\,dG^{(1)}_{p_1,p_2}(
\lambda):=p_1 \biggl(\int\phi(\lambda)\,dF^{\bbS_{{\mathbf x}{\mathbf
y}}}(\lambda)-
\int\phi(\lambda)\,dF^{c_{1n},c_{2n}}(\lambda) \biggr),
\end{equation}
and
%
%e2.14 #&#
\begin{equation}
\label{cen1} \qquad\int\phi(\lambda)\,dG^{(2)}_{p_1,p_2}(
\lambda):=p_1 \biggl(\int\phi(\lambda)\,dF^{\bbT_{{\mathbf x}{\mathbf
y}}}(\lambda)-
\int\phi(\lambda)\,dF^{c'_{1n},c'_{2n}}(\lambda) \biggr).
\end{equation}
Also, let
%
%e2.15 #&#
\begin{eqnarray}\label{b4}
\bar y_1 &:=& \frac{c_1}{1-c_2}\in(0, +\infty), \qquad
\bar y_2 := \frac{c_1}{c_2}\in(0,1),\nonumber
\\
h &=& \sqrt{\bar y_1+\bar y_2- \bar y_1\bar y_2},\qquad a_1=\frac{(1-h)^2}{(1-\bar y_2)^2},\qquad
a_2 = \frac{(1+h)^2}{(1-\bar y_2)^2},\hspace*{-20pt}
\\
g_{\bar y_1, \bar y_2}(\lambda) &=& \frac{1-\bar y_2}{2\pi\lambda(\bar y_1+\bar y_2\lambda)}\sqrt
{(a_2-\lambda) (\lambda-a_1)},\qquad
a_1<\lambda<a_2.\nonumber
\end{eqnarray}

%th2.2 #&#
\begin{teo}\label{teo1}
Let $\phi_1,\ldots,\phi_s$ be functions analytic in an open region
in the complex plane containing the interval $[a_1,a_2]$. %which is the
%support of the LSD $F_{y_1,y_2}$ defined in (\ref{FmatrixLSD}).
In addition to Assumptions \ref{ass1}~and~\ref{ass3}, suppose that
%
%e2.16 #&#
\begin{equation}
\label{b5} EW_{11}^{4}=3.
\end{equation}
%
%Under the assumption that $y_1\in(0,+\infty)$, $y_2\in(0,1)$, we
%have,
Then, as $n\rightarrow\infty$, the random vector
%
%e2.17 #&#
\begin{equation}
\biggl(\int\phi_1(\lambda)\,dG^{(1)}_{p_1,p_2}(
\lambda), \ldots, \int\phi_{s}(\lambda)\,dG^{(1)}_{p_1,p_2}(
\lambda) \biggr)
\end{equation}
converges weakly to a Gaussian vector $(X_{\phi_1},\ldots,X_{\phi
_s})$ with mean
\begin{eqnarray}
\label{mean1} EX_{\phi_i}&=&\lim_{r\downarrow1}\frac{1}{4\pi i}
\oint_{\llvert \xi
\rrvert =1}f_i \biggl(\frac{1+h^2+2h\mathfrak{R}(\xi)}{(1-\bar y_2)^2}
\biggr)
\nonumber\\[-8pt]\\[-8pt]\nonumber
&&\hspace*{62pt}{}\times  \biggl[
\frac{1}{\xi-r^{-1}}+\frac{1}{\xi+r^{-1}}-\frac{2}{\xi
+ \bar y_2/h} \biggr]\,d\xi,
\nonumber
\end{eqnarray}
and covariance function
%
%e2.18 #&#
\begin{eqnarray}\label{variance1}
&& \operatorname{cov}(X_{\phi_i},X_{\phi_j})\nonumber
\\
&&\qquad =-\lim
_{r\downarrow1}\frac{1}{4\pi
^2}\oint_{\llvert \xi_1\rrvert =1}
\oint_{\llvert \xi_2\rrvert =1}
\biggl(\biggl(f_i \biggl(\frac
{1+h^2+2h\mathfrak{R}(\xi_1)}{(1-\bar y_2)^2} \biggr)
\nonumber\\[-8pt]\\[-8pt]\nonumber
&&\hspace*{151pt}{}\times f_j \biggl(\frac
{1+h^2+2h\mathfrak{R}(\xi_2)}{(1-\bar y_2)^2} \biggr)\biggr)
\\
&&\hspace*{219pt}{}\Big/(\xi_1-r\xi_2)^2\biggr)\,d\xi_1\,d\xi_2,
\nonumber
\end{eqnarray}
where $f_i(\lambda)=\phi_i(\frac{1}{1+(({1-c_2})/{c_2})\lambda
})$, $\mathfrak{R}$ denotes the real part of a complex number and
$r\downarrow1$ means that $r$ approaches to $1$ from above.
\end{teo}

%Assumptions $1$ and $3$ can be further relaxed. For example, consider
%a general form of $\bbx=\bbT_1\bbw+\bbf$ and $\bby=\bbT_2\bbv+\bbg$,
%where $\bbT_1\dvtx  p_1\times p_1$ and $\bbT_2\dvtx  p_2\times p_2$ are
% nonsingular matrices; $\bbw\dvtx  p_1\times1$ and $\bbv\dvtx  p_2\times1$ are
%random vectors consisting of i.i.d. components with zero means and
%unit variances;
% $\bbf\dvtx  p\times1$ and $\bbg\dvtx  q\times1$ are deterministic vectors. By
%Theorem $11.1.2$ of \cite{Fujikoshi2010}, the canonical correlation
%coefficients between $\bbx$ and $\bby$ are
% the same as those between $\bbw$ and $\bbv$. % where $\bbx\dvtx  p_1\times
%1$, $\bby\dvtx  p_2\times1$, $\bbu\dvtx  p_1\times1$ and $\bbv\dvtx  p_2\times1$
%are random vectors whose sample matrices
% are $\bbX, \bbY, \bbU$ and $\bbV$.
%From this point, it is sufficient to consider the case under
%assumptions (a) and (b). In other words, the asymptotic theorems of
%canonical correlation coefficients for two random vectors with
%dependent elements are the same as those with corresponding
%independent elements.

%re3 #&#
\begin{rmk}
When $\phi(x)=x$, the mean of the limit distribution in Theorem~\ref
{teo1} is $0$ and the variance is $\frac
{2h^2y_1^2y_2^2}{(y_1+y_2)^4}$. These are calculated in Example 4.2 of
\cite{zsr}. Moreover, assumption (\ref{b5}) can be replaced by
$EY_{11}^4=3$ since $\bbX$ and $\bbY$ have an equal status in the
matrix $\bbS_{{\mathbf x}{\mathbf y}}$.
\end{rmk}

Before stating the CLT of the linear spectral statistics for the matrix
$\bbT_{{\mathbf x}{\mathbf y}}$, we make some notation.
%The Stieltjes transforms of the LSD's of $\bbT_{\bbx\bby}$
%and $\widetilde\bbP_{\bby}\bbP_{t\bbx}\bbP_{\bby}$ are denoted by $m(z)$
%and $\underline{m}(z)$ respectively.
%Set
%$$
%y_1=\frac{c_1'}{c_2'}\in(0,+\infty), y_2=
%$$
%$$
%$$
Let $r$ be a positive integer and introduce
\begin{eqnarray*}
m_r(z)&=&\int\frac{dH_t(x)}{(x-z+\varpi(z))^r},\qquad\varpi(z)=\frac
{1}{1+y_2m_{y_2t}(z)},
\\
g(z) &=& \frac{y_2(m_{y_2t}(-\underline{m}_{{\mathbf
y}}(z)))'}{(1+y_2m_{y_2t}(-\underline{m}_{{\mathbf y}}(z)))^2},
\\
s(z_1,z_2)&=&\frac{1}{1+y_2m_{y_2t}(z_1)}-\frac
{1}{1+y_2m_{y_2t}(z_2)},
\\
h(z) &=& \frac{-\underline{m}_{{\mathbf
y}}^2(z)}{1-y_1\underline{m}_{{\mathbf y}}^2(z)\int
(dF_{y_2t}(x))/{(x+\underline{m}_{{\mathbf y}}(z))^2}},
\end{eqnarray*}
where $(m_{y_2t}(z))'$ stands for the derivative with respect to $z$.

%th2.3 #&#
\begin{teo}\label{teo3}
Let $\phi_1,\ldots,\phi_s$ be functions analytic in an open region
in the complex plane containing the support of the LSD $\widetilde
F(\lambda)$ whose Stieltjes transform is (\ref{bb2}). %which is the
%support of the LSD $F_{y_1,y_2}$ defined in (\ref{FmatrixLSD}).
In addition to Assumptions \ref{ass2}--\ref{ass4}, suppose that the spectra norm of
$\bolds{\Sigma}_{{\mathbf x}{\mathbf x}}$ is bounded and
%
%e2.19 #&#
\begin{equation}
\label{2b5} EW_{11}^{4}=3.
\end{equation}

%Under the assumption that $y_1\in(0,+\infty)$, $y_2\in(0,1)$, we
%have,
%Denote
%y_{1n}=\frac{p_1}{p_2}\rightarrow y_1\in(0,+\infty), y_{2n}=

\begin{longlist}[(a)]
\item[(a)] If $c_2'\in(0,1)$, then the random vector
%
%e2.20 #&#
\begin{equation}
\label{c7} \biggl(\int\phi_1(\lambda)\,dG^{(2)}_{p_1,p_2}(
\lambda), \ldots, \int\phi_{s}(\lambda)\,dG^{(2)}_{p_1,p_2}(
\lambda) \biggr)
\end{equation}
converges weakly to a Gaussian vector $(X_{\phi_1},\ldots,X_{\phi
_s})$ with mean
%
%e2.21 #&#
\begin{eqnarray}
\label{w6} EX_{\phi_i}&=&-\frac{1}{2\pi i}\oint_{\mathcal{C}}
\phi_i \biggl(\frac{qz}{1+qz} \biggr)\nonumber
\\
&&\hspace*{39pt}{}\times \biggl(\biggl(y_1\int\underline
{m}_{{\mathbf y}}(z)^3x\bigl[x+\underline{m}_{{\mathbf y}}(z)\bigr]^{-3}\,dF_{y_2t}(x)\biggr)\nonumber
\\
&&\hspace*{60pt}{}\Big/
\biggl[1-y_1\int\underline{m}_{{\mathbf y}}(z)^2\bigl(x+\underline{m}_{{\mathbf
y}}(z)\bigr)^{-2}\,dF_{y_2t}(x)\biggr]^2
\nonumber
\\
&&\hspace*{58pt}{}+h(z) \bigl(y_2\varpi^2\bigl(-\underline{m}_{{\mathbf y}}(z)\bigr)m_3\bigl(-\underline
{m}_{{\mathbf y}}(z)\bigr)
\nonumber\\[-8pt]\\[-8pt]\nonumber
&&\hspace*{93pt}{}+y_2^2\varpi^4\bigl(-\underline{m}_{{\mathbf
y}}(z)\bigr)m_{y_2t}'\bigl(-\underline{m}_{{\mathbf y}}(z)\bigr)m_3\bigl(-\underline
{m}_{{\mathbf y}}(z)\bigr)
\bigr)
\\
&&\hspace*{91pt}{}
/ \bigl(1-y_2\varpi^2\bigl(-\underline{m}_{{\mathbf y}}(z)\bigr)m_2\bigl(-\underline{m}_{{\mathbf
y}}(z)\bigr)\bigr)
\nonumber
\\
&&\hspace*{58pt}{}-h(z) \bigl(y_2^2\varpi^3\bigl(-\underline{m}_{{\mathbf
y}}(z)\bigr)m_{y_2t}'\bigl(-\underline{m}_{{\mathbf y}}(z)\bigr)m_2\bigl(-\underline
{m}_{{\mathbf y}}(z)\bigr)\bigr)\nonumber
\\
&&\hspace*{142pt}{}/
\bigl(1-y_2\varpi^2\bigl(-\underline{m}_{{\mathbf
y}}(z)\bigr)m_2\bigl(-\underline{m}_{{\mathbf y}}(z)\bigr)\bigr) \biggr)\,dz
\nonumber
\end{eqnarray}
and covariance
\begin{eqnarray}
\label{w7}
\hspace*{-5pt}&& \operatorname{cov}(X_{\phi_i},X_{\phi_j})\nonumber
\\
\hspace*{-5pt}&&\!\qquad = -\frac{1}{2\pi^2}\nonumber
\\
\hspace*{-5pt}&&\hspace*{40pt}\times
\oint_{\mathcal
{C}_1}\oint_{\mathcal{C}_2}\phi_i \biggl(
\frac{qz_1}{1+qz_1} \biggr)\phi_j\biggl(\frac{qz_2}{1+qz_2}\biggr)
\nonumber\\[-8pt]\hspace*{-5pt} \\[-8pt]\nonumber
\hspace*{-5pt}&&\hspace*{75pt}{}\times
\biggl(\frac{\underline{m}_{{\mathbf y}}'(z_1)\underline{m}_{{\mathbf
y}}'(z_2)}{(\underline{m}_{{\mathbf y}}(z_1)-\underline{m}_{{\mathbf
y}}(z_2))^2}-\frac{1}{(z_1-z_2)^2}\nonumber
\\
\hspace*{-5pt}&&\hspace*{93pt}{}-\frac{h(z_1)h(z_2)}{(-\underline{m}_{{\mathbf y}}(z_2)+\underline
{m}_{{\mathbf y}}(z_1))^2}\nonumber
\\
\hspace*{-5pt}&&\hspace*{93pt}{}+ \frac{h(z_1)h(z_2)[1+g(z_1)
+g(z_2)+g(z_1)g(z_2)]}{[-\underline{m}_{{\mathbf y}}(z_2)+\underline
{m}_{{\mathbf y}}(z_1)+s(-\underline{m}_{{\mathbf y}}(z_1),-\underline
{m}_{{\mathbf y}}(z_2))]^2} \biggr)\,dz_1\,dz_2.
\nonumber
\end{eqnarray}
Here, $q$ is defined in Theorem \ref{teo2}. The contours in (\ref
{w6}) and (\ref{w7}) [two in (\ref{w7}), which may be assumed to be
nonoverlapping] are closed and are taken in the positive direction in
the complex plain, each enclosing the support of $\widetilde
{F}(\lambda)$.

\item[(b)] If\vspace*{1pt} $c_2'\in[1,+\infty)$ ($p_2\geq n$), (\ref{c7})
converges weakly to a Gaussian vector $(X_{\phi_1},\ldots,X_{\phi
_s})$ with mean
\begin{eqnarray}\label{cca31}
EX_{\phi_i}&=&-\frac{1}{2\pi i}
\nonumber\\[-8pt]\\[-8pt]\nonumber
&&\hspace*{6pt}{}\times \oint_{\mathcal{C}}
\phi_i \biggl(\frac
{t^{-1}z}{1+t^{-1}z} \biggr)\frac{c_1'\int(1+\lambda\underline
{s}(z)^3)^{-3}\underline{s}(z)^{3}\lambda^2\,dH(\lambda)}{
(1-c_1'\int\underline{s}(z)^{2}\lambda^2(1+\lambda\underline
{s}(z))^{-2}\,dH(\lambda) )^2}\,dz\hspace*{-25pt}
\end{eqnarray}
and
%
%e2.22 #&#
\begin{eqnarray}\label{cca32}
&& \operatorname{cov}(X_{\phi_i}, X_{\phi_j})\nonumber
\\
&&\qquad = -\frac{1}{2\pi^2}
\oint_{\mathcal
{C}_1}\oint_{\mathcal{C}_2}\phi_i \biggl(
\frac{t^{-1}z_1}{1+t^{-1}z_1} \biggr)
\\
&&\hspace*{91pt}{}\times \phi_i \biggl(\frac
{t^{-1}z_2}{1+t^{-1}z_2} \biggr)
\frac{\underline{s}{}'(z_1)\underline
{s}{}'(z_2)}{ (\underline{s}(z_1)-\underline{s}(z_2)
)^2}\,dz_1\,dz_2,\nonumber
\end{eqnarray}
where $\underline{s}(z)$ is Stieltjes transform of the LSD of the
matrix $\frac{1}{n}\bbW^{T}\bolds{\Sigma}_{{\mathbf x}{\mathbf x}}\bbW
$. The contours in (\ref{cca31}) and (\ref{cca32}) [two in (\ref
{cca32}), which may be assumed to be nonoverlapping] are closed and
are taken in the positive direction in the complex plain, each
enclosing the support of $\widetilde{G}(\lambda)$.
\end{longlist}
\end{teo}

Here, we would like to point out that the idea of testing independence
between two random vectors ${\mathbf x}$ and ${\mathbf y}$ by CCA is based on
the fact that the uncorrelatedness between ${\mathbf x}$ and ${\mathbf y}$ is
equivalent to independence between them when the random vector of size
$(p_1+p_2)$ consisting of the components of ${\mathbf x}$ and ${\mathbf y}$ is
a Gaussian random vector. See Wilks \cite{wi} and Anderson \cite{Ader1984}. For
non-Gaussian random vectors ${\mathbf x}$ and ${\mathbf y}$, uncorrelatedness
is not equivalent to independence. CCA may fail in this case. Yet,
since Theorems \ref{teo1} and \ref{teo3} hold for non-Gaussian random
vectors ${\mathbf x}$ and ${\mathbf y}$ CCA can be still utilized to capture
dependent but uncorrelated ${\mathbf x}$ and ${\mathbf y}$ such as ARCH type of
dependence by considering the high power of their entries. See
Section~\ref{secARCH} for further discussion.

%As will be seen in the Appendix, when $\bbx$ and $\bby$ are
%independent and both consist of Gaussian random variables, Theorem
%importance of lies in the fact that it can be utilized to capture
%uncorrelated but dependent $\bbx$ and $\bby$ such as ARCH type of
%dependence by considering the high power of their entries.

Following \cite{LP}, condition (\ref{b5}) can be removed. However, it
will significantly increase the length of this work and we will not
pursue it here.

%Additional to the assumptions in Proposition \ref{prop1} except that $
%and
%for any fixed $\varepsilon>0$, where $F^{(n)}_{kj}$ is the probability
%law of $X_{kj}$.

%Then as $n\rightarrow\infty$, the random vector
%converges weakly to a Gaussian vector $(X_{\phi_1},\ldots,X_{\phi_s})$
%with means (\ref{mean1}) and covariances (\ref{variance1}).

%s3 #&#
\section{Test statistics}\label{sec3}
Note that the regularized statistic $
\int\lambda \,dG^{(2)}_{p_1,p_2}(\lambda)
$ in (\ref{cen1}) [when $\phi(\lambda)=\lambda$] involves the
unknown covariance matrix $\bolds{\Sigma}_{{\mathbf x}{\mathbf x}}$
through $F^{c'_{1n},c'_{2n}}(\lambda)$. In order to apply it to
conduct tests, one needs to estimate the unknown parameter. It is well
known that estimating the population
covariance matrix $\bolds{\Sigma}_{{\mathbf x}{\mathbf x}}$ is very
challenging unless it is sparse. \cite{Karoui2008} and \cite
{Baietal2010} proposed some approaches to estimate the limit of the ESD
of $\bolds{\Sigma}_{{\mathbf x}{\mathbf x}}$ or its moments. However,
the convergence rate is not fast enough to offset the order of $p_1$.
Indeed, Theorem~1 of \cite{Baietal2010} implies that the best possible
convergence rate is $O_p(\frac{1}{n})$. In view of this, we provide
two methods to deal with the problem. One is to estimate $\int\lambda
\,dF^{c'_{1n},c'_{2n}}(\lambda)$ in a framework of sparsity while
the other one is to eliminate this unknown parameter by dividing the
samples into two groups.

%s3.1 #&#
\subsection{Plug-in estimator under sparsity}\label{sec3.1}
When $c_2'<1$, it turns out that
%
%e3.1 #&#
\begin{equation}
\label{b8*}\int\lambda \,dF^{c_{1n}',c_{2n}'}(\lambda)=\frac
{p_2}{p_1}-
\frac{p_2}{p_1}\frac{1}{1+c_{1n}'m_{nt}},
\end{equation}
where $m_{nt}$ is a solution to the equation
%
%e3.2 #&#
\begin{equation}
\label{add11y} m_{nt}=a_n-\frac{a_nt}{p_1}\tr
\bigl(a_n^{-1}\bolds{\Sigma}_{{\mathbf
x}{\mathbf x}}+t\bbI
\bigr)^{-1}
\end{equation}
with $a_n=1+c_{1n}'m_{nt}$ (see the proof of Theorem \ref{teo5}). An
estimator of $m_{nt}$ is then proposed as $\hat{m}_{nt}$ which is a
solution to the equation
%
%e3.3 #&#
\begin{equation}
\label{b9} \hat{m}_{nt}=\hat{a}_n-
\frac{\hat{a}_nt}{p_1}\tr \bigl(\hat{a}_n^{-1}\hat{\bolds{
\Sigma}}_{{\mathbf x}{\mathbf x}}+t\bbI\bigr)^{-1},
\end{equation}
with $\hat{a}_n=1+c_{1n}'\hat{m}_{nt}$.
Here, we use a thresholding estimator $\hat{\bolds{\Sigma
}}_{{\mathbf x}{\mathbf x}}$ to estimate $\bolds{\Sigma}_{{\mathbf x}{\mathbf
x}}$, slightly different from that proposed by
\cite{BL20081}. Specifically speaking, suppose that
the underlying random variables $\{X_{ij}\}$ are mean zero and variable
one. Then define $\hat{\bolds{\Sigma}}_{{\mathbf x}{\mathbf x}}$ to be
a matrix whose diagonal entries are all one and the off diagonal
entries are $\hat{\sigma}_{ij}I(\llvert \hat{\sigma}_{ij}\rrvert \geq
\ell)$
with $\ell=M\sqrt{\frac{\log p_1}{n}}$ and $M$ being some
appropriate constant ($M$ will be selected by cross-validation). Here
$\hat{\sigma}_{ij}$ denotes the entry at the $(i,j)$th position of
sample covariance matrix $\frac{1}{n}\bbX\bbX^T$. Therefore, the
resulting test statistic is
%
%e3.4 #&#
\begin{equation}
\label{b16} p_1 \biggl(\int\lambda \,dF^{\bbT_{{\mathbf x}{\mathbf y}}}(\lambda)-
\biggl(\frac
{p_2}{p_1}-\frac{p_2}{p_1}\frac{1}{1+c_{1n}'\hat{m}_{nt}}\biggr) \biggr).
\end{equation}

When $p_2\geq n$, it turns out that
%
%e3.5 #&#
\begin{equation}
\label{add75y} \int\lambda \,dF^{c_{1n}',c_{2n}'}(\lambda)=1-tm_n^{(1t)},
\end{equation}
where $m_n^{(1t)}$ satisfies the equation
%
%e3.6 #&#
\begin{equation}
\label{add10y} m_n^{(1t)}=\frac{1}{p_1}\tr  \bigl(
\bigl(1-c_{1n}'+c_{1n}'tm_n^{(1t)}
\bigr)\bolds{\Sigma}_{{\mathbf x}{\mathbf
x}}+t\bbI\bigr)^{-1}.
\end{equation}
We then propose the resulting test statistic
%
%e3.7 #&#
\begin{equation}
\label{b17} p_1 \biggl(\int\lambda \,dF^{\bbT_{{\mathbf x}{\mathbf y}}}(\lambda)-
\bigl(1-t\hat{m}_n^{(1t)}\bigr) \biggr),
\end{equation}
where $\hat{m}_n^{(1t)}$ satisfies the equation
%
%e3.8 #&#
\begin{equation}
\label{b15*} \hat{m}_n^{(1t)}=\frac{1}{p_1}\tr  \bigl(
\bigl(1-c_{1n}'+c_{1n}'t\hat
{m}_n^{(1t)}\bigr)\hat{\bolds{\Sigma}}_{{\mathbf x}{\mathbf x}}+t
\bbI\bigr)^{-1}.
\end{equation}

%th3.1 #&#
\begin{teo}\label{teo5}
In addition to assumptions in Theorem \ref{teo3}, suppose that
$EX_{ij}^2=1$, $\sup _{i,j}E\llvert X_{ij}\rrvert ^{17}<\infty$
for all $i$
and $j$ and that
%
%e3.9 #&#
\begin{equation}
\label{b18} s_o(p_1) \biggl(\frac{\log p_1}{n}
\biggr)^{(1-q)/2}\rightarrow0,
\end{equation}
where $
\sum _{i\neq j}\llvert \sigma_{ij}\rrvert ^q=s_o(p_1)$ with
$0\leq q<1$.
\begin{longlist}[(a)]
\item[(a)] If $c_2'<1$, then $p_1 (\int\lambda \,dF^{\bbT_{{\mathbf x}{\mathbf
y}}}(\lambda)-(\frac{p_2}{p_1}-\frac{p_2}{p_1}\frac
{1}{1+c_{1n}'\hat{m}_{nt}}) )$ converges weakly to a normal
distribution with the mean and variance given in (\ref{w6}) and (\ref
{w7}) with $\phi(\lambda)=\lambda$.

\item[(b)] If $c_2'\geq1$, then $p_1 (\int\lambda \,dF^{\bbT_{{\mathbf x}{\mathbf
y}}}(\lambda)-(1-t\hat{m}_n^{(1t)}) )$ converges weakly to a
normal distribution with the mean and variance given in part \textup{(b)} of
Theorem \ref{teo3} with $\phi(\lambda)=\lambda$.
\end{longlist}
\end{teo}
We demonstrate an example of sparse covariance matrices in the
simulation parts, satisfying the sparse condition (\ref{b18}).
%Instead of assuming $EX_{ij}^2=1$ we may also suppose that $EX_{ij}^2=
%be $\bolds{\Sigma}_{\bbx\bbx}/\sigma$. Note that $\bbX^T(\bbX

%s3.2 #&#
\subsection{Strategy of dividing samples}\label{sec3.2}
If (\ref{b18}) is not satisfied, we then propose a strategy of
dividing the total samples into two groups. Specifically speaking, we
divide the $n$ samples of $({\mathbf x}, {\mathbf y})$ into two groups,
respectively, that is,
%
%e3.10 #&#
\begin{equation}
\label{add50y}
\mbox{\textit{Group} 1:}\quad \bbX^{(1)}=({\mathbf x}_1, {\mathbf
x}_2, \ldots, {\mathbf x}_{[n/2]}),\qquad \bbY^{(1)}=({\mathbf
y}_1, {\mathbf y}_2, \ldots, {\mathbf y}_{[n/2]})\hspace*{-20pt}
\end{equation}
and
%
%e3.11 #&#
\begin{eqnarray}
\label{add51y}
\mbox{\textit{Group} 2:}\quad \bbX^{(2)}&=&({\mathbf x}_{[n/2]+1}, {\mathbf
x}_{[n/2]+2}, \ldots, {\mathbf x}_{n}),
\nonumber\\[-8pt]\\[-8pt]\nonumber
\bbY^{(2)}&=&({\mathbf
y}_{[n/2]+1}, {\mathbf y}_{[n/2]+2}, \ldots, {\mathbf y}_n),
\end{eqnarray}
where $[n/2]$ is the largest integer not greater than $n/2$. When $n$
is odd, we discard the last sample. However, if the above strategy of
dividing samples into two groups is directly used, then the asymptotic
means of the resulting statistic [the difference between the statistics
in (\ref{cen1}) obtained from two subsamples] are always zero in both
null hypothesis and alternative hypothesis due to similarity of two
groups so that the power of the test statistic is very low. This is
also confirmed by simulations. Therefore, we further propose its
modified version as follows.

For $\bbY^{(2)}$ in Group 2, we extract a sub-data $\widetilde{\bbY
}{}^{(2)}$, that is,
\[
\widetilde{\bbY}{}^{(2)}=(\tilde{{\mathbf y}}_{[n/2]+1}, \tilde
{{\mathbf y}}_{[n/2]+2}, \ldots, \tilde{{\mathbf y}}_n),
\]
where $\tilde{{\mathbf y}}_j$ consists of the first $[p_2/2]$
components of ${\mathbf y}_j$, for all $j=[n/2]+1$, $[n/2]+2, \ldots, n$. We
use $\widetilde{\bbY}{}^{(2)}$ to form a new group
\begin{eqnarray*}
\mbox{\textit{Modified Group} 2:}\quad \bbX^{(2)} &=& ({\mathbf x}_{[n/2]+1}, {\mathbf
x}_{[n/2]+2}, \ldots, {\mathbf x}_{n}),
\\
\widetilde{\bbY }{}^{(2)}&=&(\tilde{{\mathbf y}}_{[n/2]+1}, \tilde{{\mathbf y}}_{[n/2]+2}, \ldots, \tilde{{\mathbf y}}_n).
\end{eqnarray*}
%

%In order to derive bigger powers under alternative hypothesis, we
%propose Modified Group 2 instead of using Group 2 \,directly. We will
%emphasize this point after introducing the modified statistics with
%the sub-sample data.

For Group 1, it follows from Theorem~\ref{teo3} that
%
%e3.12 #&#
\begin{equation}
\label{add1y} \int\lambda \,dp_1 \bigl(F^{\bbT_{{\mathbf x}{\mathbf
y}}^{(1)}}(\lambda
)-F^{2c_{1n}', 2c_{2n}'}(\lambda) \bigr)\stackrel{d} {\rightarrow}Z_1,
\end{equation}
where $\bbT_{{\mathbf x}{\mathbf y}}^{(1)}$ is obtained from $\bbT_{{\mathbf
x}{\mathbf y}}$ with $\bbX$ and $\bbY$ replaced by $\bbX^{(1)}$ and
$\bbY^{(1)}$, respectively, and $Z_1$ is a normal random variable with
mean and variance given in Theorem \ref{teo3} with $c_1'$ and $c_2'$
replaced by $2c_1'$ and $2c_2'$, respectively, and $\phi(\lambda
)=\lambda$. Similarly, with Modified Group 2, by Theorem~\ref{teo3}
%
%e3.13 #&#
\begin{equation}
\label{add2y} \int\lambda \,dp_1 \bigl(F^{\bbT_{{\mathbf x}{\mathbf
y}}^{(2)}}(\lambda
)-F^{2c_{1n}', c_{2n}'}(\lambda) \bigr)\stackrel{d} {\rightarrow}Z_2,
\end{equation}
where $\bbT_{{\mathbf x}{\mathbf y}}^{(2)}$ is $\bbT_{{\mathbf x}{\mathbf y}}$ with
$\bbX$ and $\bbY$ replaced by $\bbX^{(2)}$ and $\widetilde{\bbY
}{}^{(2)}$, respectively, and $Z_2$ is a normal random variable with the
mean and variance given in Theorem \ref{teo3} with $\phi(\lambda
)=\lambda$ and $c_1'$ replaced by $2c_1'$.

We next investigate the relation between
\[
\int\lambda \,dF^{2c_{1n}', 2c_{2n}'}(\lambda)\quad\mbox{and}\quad \int\lambda \,dF^{2c_{1n}',
c_{2n}'}(\lambda),
\]
and then calculate some difference between the
two statistics in~(\ref{add1y}) and (\ref{add2y}) in~order to
eliminate the unknown parameters $\int\lambda \,dF^{2c_{1n}',
2c_{2n}'}(\lambda)$ and $\int\lambda \,dF^{2c_{1n}',
c_{2n}'}(\lambda)$.

When $c_2'<1/2$, we have
%
%e3.14 #&#
\begin{equation}
\label{add3y} \int\lambda \,dF^{2c_{1n}',2c_{2n}'}(\lambda)=\frac
{p_2}{p_1}-
\frac{p_2}{p_1}\frac{1}{1+2c_{1n}'\tilde{m}_{nt}},
\end{equation}
where $\tilde{m}_{nt}$ is obtained from $m_{nt}$ satisfying (\ref
{add11y}) with $c_{1n}'$ replaced by $2c_{1n}'$. On the other hand,
%
%e3.15 #&#
\begin{equation}
\label{add3y*} \int\lambda \,dF^{2c_{1n}',c_{2n}'}(\lambda)=\frac
{p_2/2}{p_1}-
\frac{p_2/2}{p_1}\frac{1}{1+2c_{1n}'\tilde{m}_{nt}}.
\end{equation}
It follows that
%
%e3.16 #&#
\begin{equation}
\label{add5y} \int\lambda \,dF^{2c_{1n}',2c_{2n}'}(\lambda)=2\int
\lambda
\,dF^{2c_{1n}',c_{2n}'}(\lambda).
\end{equation}

When $[p_2/2]>[n/2]$, we have
%
%e3.17 #&#
\begin{equation}
\label{add12y} \int\lambda \,dF^{2c_{1n}',2c_{2n}'}(\lambda)=\int
\lambda
\,dF^{2c_{1n}',c_{2n}'}(\lambda)=1-t\tilde{m}_n^{(1t)},
\end{equation}
where $\tilde{m}_n^{(1t)}$ is $m_n^{(1t)}$ satisfying (\ref
{add10y}) with $c_{1n}'$ replaced by $2c_{1n}'$.

The last case is $[p_2/2]\leq[n/2]$ and $c_2'\geq1/2$. For this case,
if we still consider Group 1 and Modified Group 2, then
\begin{eqnarray*}
\int\lambda \,dF^{2c_{1n}',2c_{2n}'}(\lambda)&=&1-t\tilde{m}_n^{(1t)},
\\
\int\lambda \,dF^{2c_{1n}',c_{2n}'}(\lambda)&=&\frac
{[p_2/2]}{p_1}-\frac{[p_2/2]}{p_1}
\frac{1}{1+2c_{1n}'\tilde{m}_{nt}}.
\end{eqnarray*}
From the above formulas, it is difficult to figure out the relation
between $\int\lambda \,dF^{2c_{1n}', 2c_{2n}'}(\lambda)$ and
$\int\lambda \,dF^{2c_{1n}', c_{2n}'}(\lambda)$ depending on the
unknown parameter $\bolds{\Sigma}_{{\mathbf x}{\mathbf x}}$. To overcome
this difficulty, we also apply a ``sub-data'' trick to Group 1.
Specifically speaking, consider a Modified Group 1 as follows.
\[
\mbox{\textit{Modified Group} 1:}\quad \bbX^{(1)}=({\mathbf x}_1, {\mathbf
x}_2, \ldots, {\mathbf x}_{[n/2]}),\qquad \dot{\bbY}^{(1)}=(
\dot{\mathbf y}_1, \dot{\mathbf y}_2, \ldots, \dot{\mathbf
y}_{[n/2]}),
\]
where $\dot{\mathbf y}_k$ consists of the last $[p_2/2]$ components of
${\mathbf y}_k$, that is, the $i$th component of $\dot{\mathbf y}_k$ is the
$([p_2/2]+i)$th component of ${\mathbf y}_k$, for all $i=1,2,\ldots,[p_2/2]$ and $k=1,2,\ldots,[n/2]$. For Modified Group 1, by Theorem~\ref{teo3}, we have
%
%e3.18 #&#
\begin{equation}
\label{add90y} \int\lambda \,dp_1 \bigl(F^{\widetilde{\bbT}_{{\mathbf x}{\mathbf
y}}^{(1)}}(
\lambda)-F^{2c_{1n}', c_{2n}'}(\lambda) \bigr)\stackrel{d}
{\rightarrow}
Z_3,
\end{equation}
where $\widetilde{\bbT}_{{\mathbf x}{\mathbf y}}^{(1)}$ is $\bbT_{{\mathbf
x}{\mathbf y}}$ with $\bbX$ and $\bbY$ replaced by $\bbX^{(1)}$ and
$\dot{\bbY}^{(1)}$, respectively; and $Z_3$ is a normal random
variable with the mean and variance given in Theorem \ref{teo3} with
$\phi(\lambda)=\lambda$ and $c_{1}'$ replaced by $2c_{1}'$.
Since the unknown parameters in (\ref{add2y}) and (\ref{add90y}) are
the same the difference between (\ref{add2y}) and (\ref{add90y}) can
be taken as the modified statistic.

The asymptotic distributions of the three resulting statistics are
given in Theorem \ref{teo6}.

%%%%%%%%%%%%%%%%%%%%%%%%%%%%%%%%%%%%%%%%%%%%%%%%%%

%th3.2 #&#
\begin{teo}\label{teo6}
%In addition to Assumptions $2$-$4$, suppose that $EX_{11}^{4}=3$.
Suppose that assumptions in Theorem \ref{teo3} hold.
\begin{longlist}[(a)]
\item[(a)] If $c_{2}'<1/2$, the statistic $\int\lambda \,dF^{\bbT_{{\mathbf
x}{\mathbf y}}^{(1)}}(\lambda)-2\int\lambda \,dF^{\bbT_{{\mathbf x}{\mathbf
y}}^{(2)}}(\lambda)$ converges\break weakly to a normal distribution with
the mean $(\mu_1-2\mu_2)$ and variance $(\sigma_1^2+4\sigma_2^2)$,
where $\mu_1$ and $\sigma_1^2$ are given in (\ref{w6}) and (\ref
{w7}), respectively, with $c_{1}', c_{2}'$ replaced by
$2c_{1}', 2c_{2}'$, respectively, and $\phi(\lambda)=\lambda$;
$\mu_2$ and $\sigma_2^2$ are given in (\ref{w6}) and (\ref{w7}),
respectively, with $c_{1}'$ replaced by $2c_{1}'$ and $\phi
(\lambda)=\lambda$.

\item[(b)] If $c_{2}'\geq1$, the statistic $\int\lambda \,dF^{\bbT_{{\mathbf
x}{\mathbf y}}^{(1)}}(\lambda)-\int\lambda \,dF^{\bbT_{{\mathbf x}{\mathbf
y}}^{(2)}}(\lambda)$ converges weakly to a normal distribution with
the mean zero and variance $2\sigma_3^2$, where $\sigma_3^2$ is given
in~(\ref{cca32}) with $c_{1}'$ replaced by $2c_{1}'$ and $\phi
(\lambda)=\lambda$.

\item[(c)] If $1/2\leq c_{2}'<1$, the statistic $\int\lambda \,dF^{\widetilde
{\bbT}_{{\mathbf x}{\mathbf y}}^{(1)}}(\lambda)-\int\lambda \,dF^{\bbT_{{\mathbf
x}{\mathbf y}}^{(2)}}(\lambda)$ converges weakly to a normal distribution
with mean zero and variance $2\sigma_4^2$, where $\sigma_4^2$ is
given in (\ref{w7}) with $c_1'$ replaced by $2c_1'$ and $\phi
(\lambda)=\lambda$.
\end{longlist}
\end{teo}

%re4 #&#
\begin{rmk}
Unlike using Group 2 of (\ref{add51y}) although the asymptotic means
of the statistics in the cases (b) and (c) are zero under the null
hypothesis, they are not necessarily equal to zero under the
alternative hypothesis so that the power of the resulting test
statistic becomes much better. % A natural question may be raised for
%the dividing-sample method: why the modified Group 2 in (\ref{add52y})
%is used instead of Group 2 in (\ref{add51y})? The reason relies on
%guaranteeing relatively big powers under the alternative hypothesis.
%If Group 2 is adopted, under any alternative hypothesis, the mean of
%the statistic $\int\lambda \,dF^{\bbT_{\bbx\bby}^{(1)}}(\lambda)-2\int
%equivalent to that under the null hypothesis. Due to this phenomenon,
%the powers will be rather low under the hypothesis.
\end{rmk}

%%%%%%%%%%%%%%%%%%%%%%%%%%%%%%%%%%%%%%%%%%%%%%%%

Theorem \ref{teo6} proposes test statistics which do not involve the
unknown parameter $H$, that is, the LSD of the matrix $\bolds
{\Sigma}_{{\mathbf x}{\mathbf x}}$. However, their asymptotic means and
asymptotic variances contain some terms involving the unknown parameter
$H$. We below provide consistent estimators for such terms appearing in
(\ref{w6})--(\ref{cca32}) in\vadjust{\goodbreak} Theorem \ref{teo3}, which, together
with Slutsky's theorem and the dominant convergence theorem, are enough
for applications.

We use an estimator developed by \cite{Karoui2008} for $H$. For easy
reference, we briefly state his estimator $\hat{H}_{p_1}$ for $H$ in
the following proposition.

%pr2 #&#
\begin{prop}[(Theorem 2 of \cite{Karoui2008})]\label{prop22}
In addition to Assumptions \ref{ass2}--\ref{ass4}, suppose that the spectra norm of
$\bolds{\Sigma}_{{\mathbf x}{\mathbf x}}$ is bounded. Let $J_1, J_2,
\ldots$ be a sequence of integers tending to $\infty$. Let $z_0\in
\mathbb{C}^{+}$ and $r\in\mathbb{R}^{+}$ be such that $\bbB
(z_0,r)\subset\mathbb{C}^{+}$, where $\bbB(z_0,r)$ denotes the
closed ball of center $z_0$ and radius $r$. Let $z_1, z_2, \ldots$ be
a sequence of complex variables with a limit point, all contained in
$\bbB(z_0,r)$. Let $\hat{H}_{p_1}$ be the solution of
%
%e3.19 #&#
\begin{equation}
\label{06021} \hat{H}_{p_1}=\arg\min_{G}\max
_{j\leq J_n}\biggl\llvert\frac{1}{\hat
{\underline{s}}_{n}(z_j)}+z_j-
\frac{p_1}{n}\int\frac{\lambda
\,dG(\lambda)}{1+\lambda\hat{\underline{s}}_{n}(z_j)}\biggr\rrvert,
\end{equation}
where $G$ is a probability measure and $\hat{\underline{s}}_{n}(z)$
is the Stieltjes transform of the ESD of the matrix $\underline{\bbA
}_{{\mathbf x}{\mathbf x}}=\frac{1}{n}\bbX^{T}\bbX$. Then we have
\[
\hat{H}_{p_1}\Rightarrow H\qquad\mbox{a.s.}
\]
\end{prop}

%re5 #&#
\begin{rmk}
The estimator $\hat{H}_{p_1}$ given in (\ref{06021}) is proposed
based on the Mar\v{c}enko--Pastur equation which links the
Stieltjes transform of the empirical spectral distribution of the
sample covariance matrix to an integral against the population spectral
distribution, that is, the LSD $H$ satisfies the Mar\v{c}enko--Pastur equation
%
%e3.20 #&#
\begin{equation}
-\frac{1}{\underline{s}(z)}=z-c\int\frac{\lambda \,dH(\lambda
)}{1+\lambda\underline{s}(z)}\qquad \forall z\in\mathbb{C}^{+},
\end{equation}
where $\underline{s}(z)$ is the limit of $\hat{\underline
{s}}_{n}(z)=\frac{1}{n}\tr (\underline{\bbA}_{{\mathbf x}{\mathbf x}}-z\bbI
_{n})^{-1}$.
\end{rmk}

El Karoui \cite{Karoui2008} developed an algorithm for the estimator $\hat{H}_{p_1}$
in (\ref{06021}) and we state it below.

\begin{longlist}[(2)]
\item[(1)] \textit{A} ``\textit{basis pursuit}'' \textit{for measure space}. Instead of
searching among all possible probability measures, the search space is
restricted to mixtures of certain classes of probability measures in
order to deal with (\ref{06021}). In other words, first select a
``dictionary'' of probability measures on the real line and then
decompose the estimator on this dictionary, searching for the best
coefficients. Hence, the problem can be formulated as
\begin{eqnarray*}
&&\mbox{{finding the best possible weights}}
\\
&&\qquad \{\hat{w}_1, \ldots, \hat{w}_K\}\qquad\mbox{with } d\hat
{H}_{p_1}=\sum^{K}_{i=1}
\hat{w}_i\,dM_i,
\sum^{K}_{i=1} \hat{w}_i=1, w_i\geq0,
\end{eqnarray*}
where the $M_i$'s are the measures in the dictionary. A ``probability
measures'' dictionary is given as follows:
\begin{enumerate}[2.]
\item[1.] Point masses $\{\delta_{t_k}(x)\}_{k=1}^{K}$, where $\{t_k\}
_{k=1}^{K}$ is a grid of points.
\item[2.] Probability measures that are uniform on an interval: in this
case,\break $dM^{(1)}_i(x)=I_{[a_i, b_i]}(x)\,dx/(b_i-a_i)$.
\item[3.] Probability measures that have a linearly increasing (or
decreasing) density on an interval $[d_i, h_i]$ and density $0$
elsewhere. So, for the increasing case, $dM^{(2)}_i(x)=I_{[d_i,
h_i]}(x)\cdot2(x-d_i)/((h_i-d_i)^2)\,dx$, and density $0$ elsewhere.
\end{enumerate}

\item[(2)] \textit{A convex optimization problem}.
Let
\[
e_j=\frac{1}{\hat{\underline{s}}_{n}(z_j)}+z_j-\frac{p_1}{n}\int
\frac{\lambda \,dG(\lambda)}{1+\lambda\hat{\underline{s}}_{n}(z_j)},\qquad
j=1,2,\ldots,J_n,
\]
where $G(\cdot)$ has the form of
\[
G(x)=\sum^{K_1}_{i=1}w_i
\delta_{t_i}(x)+\sum^{K_2}_{i=K_1+1}w_i\,dM^{(1)}_i(x)+
\sum^{K}_{i=K_2+1}w_i\,dM^{(2)}_i(x),
\]
with all points $t_k, k=1,2,\ldots,K_1$, intervals $[a_i, b_i],
i=K_1+1, \ldots, K_2$ and intervals $[d_j, h_j], j=K_2+1, \ldots, K$
being in the interval $[\ell_{p_1}, \ell_1]$. Here, $\ell_{p_1}$ and
$\ell_1$ are, respectively, the smallest and largest eigenvalues of
the sample covariance matrix $\bbA_{{\mathbf x}{\mathbf x}}=\frac{1}{n}\bbX
\bbX^{T}$. Moreover, $1<K_1<K_2<K$.

The ``translation'' of the problem (\ref{06021}) into a convex
optimization problem is
\begin{eqnarray*}
&\displaystyle \min_{w_1,\ldots,w_K,u}u&
\\
&\displaystyle \forall j=1,\ldots,J_n,\qquad -u\leq\Re(e_j)\leq u,&
\\
&\displaystyle \forall j=1,\ldots,J_n,\qquad -u\leq\Im(e_j)\leq u &
\\
&\displaystyle \mbox{subject to}\qquad\sum^{K}_{i=1}w_i=1\mbox{ and }w_i\geq0, \forall i=1,2,\ldots,K.&
\end{eqnarray*}

The following proposition provides consistency of the proposed
algorithm above.
\end{longlist}

%
%pr3 #&#
\begin{prop}[(Corollary 1 of \cite{Karoui2008})]\label{prop33}
Assume the same assumptions as in Proposition \ref{prop22}. Call $\hat
{H}_{p_1}$ the solution of (\ref{06021}), where the optimization is
now over measures which are sums of atoms, the locations of which are
restricted to belong to a grid (depending on $n$) whose step size is
going to $0$ as $n\rightarrow\infty$. Then
\[
\hat{H}_{p_1}\Rightarrow H\qquad\mbox{a.s.}
\]
\end{prop}

%%%%%%%%%%%%%%%%%%%%%%%%%%%%%%%%%%%%%%%%%%%%%%%%%
Then the estimator $\hat{H}_{p_1}$ derived from the algorithm has the
form of
%
%e3.21 #&#
\begin{equation}\label{06022}
\qquad d\hat{H}_{p_1}(x)=\sum^{K_1}_{i=1}
\hat{w}_i\delta_{t_i}(x)+\sum^{K_2}_{i=K_1+1}
\hat{w}_i\,dM^{(1)}_i(x)+\sum
^{K}_{i=K_2+1}\hat{w}_i\,dM^{(2)}_i(x),
\end{equation}
where $\{\hat{w}_i, i=1,2,\ldots,K\}$ is the solution of the convex
optimization problem above. In practice, we follow the implementation
details in Appendix of \cite{Karoui2008} to derive the estimator $\hat
{H}_{p_1}(x)$.

Once the estimator (\ref{06022}) of the LSD $H$ is available, we are
in a position to provide estimators for $m_{H}(z)$ and $m_{y_2t}(z)$,
the Stieltjes transforms of $H$ and the LSD of the matrix $\bbS_{2t}$,
respectively.

%pr4 #&#
\begin{prop}\label{prop2}\label{pr4}
For any $z\in\mathcal{C}_1$, the contour specified in Theorem~\ref{teo3}, let
%
%e3.22 #&#
\begin{eqnarray}
\label{06031} \hat{m}_{H}(z)&=&\sum^{K_1}_{i=1}
\frac{\hat{w}_i}{t_i-z}+\sum^{K_2}_{i=K_1+1}
\frac{\hat{w}_i}{b_i-a_i}\log\frac
{b_i-z}{a_i-z}
\nonumber\\[-8pt]\\[-8pt]\nonumber
&&{}+\sum^{K}_{i=K_2+1}\frac{2\hat{w}_i}{(h_i-d_i)^2}
\biggl(h_i-d_i+(z-d_i)\log
\frac{h_i-z}{d_i-z} \biggr),
\end{eqnarray}
where $\log$ stands for the corresponding principal branch. The
estimator $\hat{m}_{y_2t}(z)$ satisfies the equation
%
%e3.23 #&#
\begin{equation}
\label{06032} \hat{m}_{y_2t}(z)=\hat{m}_{H_t} \biggl(z-
\frac{1}{1+y_{2n}\hat
{m}_{y_2t}(z)} \biggr),
\end{equation}
where $\hat{m}_{H_t}(z)$ is
%
%e3.24 #&#
\begin{equation}
\label{06033} \hat{m}_{H_t}(z)=-\frac{1}{z}-\frac{t}{(1-c_{2n}')z^2}
\hat{m}_{H} \biggl(\frac{t}{(1-c_{2n}')z} \biggr).
\end{equation}
Then $\hat{m}_{H}(z)$ in (\ref{06031}) and $\hat{m}_{y_2t}(z)$ in
(\ref{06032}) are consistent estimators of $m_{H}(z)$ and
$m_{y_2t}(z)$, respectively.
\end{prop}

%re6 #&#
\begin{rmk}
As to choices of intervals $[a_i, b_i]\dvtx  i=K_1+1, K_1+2, \ldots, K_2$;
and $[h_j, d_j]\dvtx  j=K_2+1, K_2+2, \ldots, K$ in (\ref{06022}), we
follow the implementation details provided in Appendix of \cite
{Karoui2008}. Furthermore, choice of $ (z_j, \hat{\underline
{s}}_n(z_j) )$ in the convex optimization problem, choice of
interval to estimate $H(\cdot)$, and choice of dictionary are all
provided in the Appendix of \cite{Karoui2008}.
\end{rmk}
%
%%%%%%%%%%%%%%%%%%%%%%%%%%%%%%%%%%%%%%%%%%%%%%%%
With consistent estimators for $m_{H}(z)$ and $m_{y_2t}(z)$ in
Proposition \ref{prop2}, we may further provide consistent estimators
for all terms appearing in (\ref{w6})--(\ref{cca32}).

%pr5 #&#
\begin{prop}\label{prop3}\label{pr5}
1. The estimator for $m_{{\mathbf y}}(z)$ is
\[
\hat{m}_{{\mathbf y}}(z)=\frac{q_n}{(1+q_nz)^2} \biggl(m_{\bbT_{{\mathbf
x}{\mathbf y}}} \biggl(
\frac{q_nz}{1+q_nz} \biggr)-1-q_nz \biggr),
\]
where $y_{1n}=\frac{c_{1n}'}{c_{2n}'}$, $y_{2n}=\frac
{c_{1n}'}{1-c_{2n}'}$, $q_n=\frac{c_{2n}'}{1-c_{2n}'}$ and
$m_{\bbT_{{\mathbf x}{\mathbf y}}}(z)$ is the Stieltjes transform of the
matrix $\bbT_{{\mathbf x}{\mathbf y}}$.

2.
The estimators for $\varpi(z)$ and $s(z_1, z_2)$ are, respectively,
\begin{eqnarray*}
\hat\varpi(z)&=&\frac{1}{1+y_{2n}\hat{m}_{y_2t}(z)},
\\
\hat{s}(z_1, z_2)&=&
\frac{1}{1+y_{2n}\hat{m}_{y_2t}(z_1)}-\frac
{1}{1+y_{2n}\hat{m}_{y_2t}(z_2)}.
\end{eqnarray*}

3. The estimators for $m_r(z)$ with $r=2,3$ are
\[
\hat{m}_2(z)=\hat{m}_{H_t}^{(1)} \bigl(z-\hat
\varpi(z) \bigr),\qquad \hat{m}_3(z)=\tfrac{1}{2}\hat{m}_{H_t}^{(2)}
\bigl(z-\hat\varpi(z) \bigr),
\]
respectively,
where $\hat{m}_{H_t}^{(j)}(z)$ is the $j$th derivative of $\hat
{m}_{H_t}(z)$ with respect to $z$ with $j=1,2$.

4.
The estimators for $g(z)$ and $h(z)$ are
\[
\hat{g}(z)=\frac{y_{2n} (\hat{m}^{(1)}_{y_2t}(-\hat{\underline
{m}}_{{\mathbf y}}(z)) )}{ (1+y_{2n}\hat{m}_{y_2t}(-\hat
{\underline{m}}_{{\mathbf y}}(z)) )^2},\qquad \hat{h}(z)=\frac{-\hat{\underline
{m}}_{{\mathbf y}}^2(z)}{1-y_{1n}\hat
{\underline{m}}_{{\mathbf y}}^2(z)\hat{m}_{y_2t}^{(1)} (-\hat
{\underline{m}}_{{\mathbf y}}(z) )},
\]
respectively,
where $\hat{\underline{m}}_{{\mathbf y}}(z)=-\frac
{1-y_{1n}}{z}+y_{1n}\hat{m}_{{\mathbf y}}(z)$.

5.~The estimators for $\varpi_1=\int\underline{m}{}^3_{{\mathbf
y}}(z)x[x+\underline{m}_{{\mathbf y}}(z)]^{-3}\,dF_{y_2t}(x)$ and $\varpi
_2=\int\underline{m}{}^2_{{\mathbf y}}(z)[x+\underline{m}_{{\mathbf
y}}(z)]^{-2}\,dF_{y_2t}(x)$ are
\begin{eqnarray*}
\hat{\varpi}_1&=&\hat{\underline{m}}^{3}_{{\mathbf y}}(z)
\hat{m}_{y_2t}^{(1)} \bigl(-\hat{\underline{m}}_{{\mathbf y}}(z)
\bigr)-\tfrac
{1}{2}\hat{\underline{m}}^{4}_{{\mathbf y}}(z)
\hat{m}_{y_2t}^{(2)} \bigl(-\hat{\underline{m}}_{{\mathbf y}}(z)
\bigr),
\\
\hat\varpi_2&=&\hat{\underline{m}}_{{\mathbf y}}^2(z)
\hat{m}_{y_2t}^{(1)} \bigl(-\hat{\underline{m}}_{{\mathbf y}}(z)
\bigr),
\end{eqnarray*}
respectively.

6.~The estimators for $\varpi_3=\int[1+\lambda\underline
{s}{}^3(z)]^{-3}\underline{s}(z)^3\lambda^2\,dH(\lambda)$ and $\varpi
_4=\int\underline{s}{}^2(z)\*\lambda^2[1+\lambda\underline
{s}(z)]^{-2}\,dH(\lambda)$ are
\begin{eqnarray*}
\hat{\varpi}_3&=&\frac{1}{[\hat{\underline{s}}_n(z)]^6}\hat{m}_{H}
\biggl(-
\frac{1}{[\hat{\underline{s}}_n(z)]^{3}} \biggr)-\frac
{2}{[\hat{\underline{s}}_n(z)]^{9}}\hat{m}_{H}^{(1)}
\biggl(-\frac
{1}{\hat{\underline{s}}_n(z)} \biggr)
\\
&&{} +\frac{1}{2[\hat{\underline
{s}}_n(z)]^{12}}\hat{m}_{H}^{(2)}
\biggl(-\frac{1}{\hat{\underline{s}}_n(z)} \biggr),
\\
\hat{\varpi}_4&=&1-\frac{2}{\hat{\underline{s}}_n(z)}\hat{m}_H
\biggl(-\frac{1}{\hat{\underline{s}}_n(z)} \biggr) +\frac{1}{[\hat
{\underline{s}}_n(z)]^2}\hat{m}_{H}^{(1)}
\biggl(-\frac{1}{\hat{\underline{s}}_n(z)} \biggr),
\end{eqnarray*}
respectively, where $\hat{\underline{s}}_n(z)$ is defined as
$\hat{\underline{s}}_n(z)=\frac{1}{n}\tr  (\underline{\bbA
}_{{\mathbf x}{\mathbf x}}-z\bbI_n )^{-1}$ with $\underline{\bbA}_{{\mathbf
x}{\mathbf x}}=\frac{1}{n}\bbX^{T}\bbX$.

All estimators listed above are consistent for the corresponding
unknown parameters.
\end{prop}

The proofs of Propositions \ref{prop2}~and~\ref{prop3} are
provided in Appendix A of \cite{supp}.

%%%%%%%%%%%%%%%%%%%%%%%%%%%%%%%%%%%%%%%%%%%%%%%%%

%s4 #&#
\section{The power under local alternatives}\label{sec4}

%Although it is difficult to provide a central limit theorem for the
%statistic $S_n$ or $T_n$ under general alternative hypothesis $
This section is to evaluate the power of $S_n$ or $T_n$ under a kind of
local alternatives. Consider the alternative hypothesis
\[
\mathbb{H}_1\dvtx  {\mathbf x}\mbox{ and }{\mathbf y}\mbox{ are dependent},
\]
satisfying condition (\ref{w4}) below. Draw $n$ samples from such
alternatives ${\mathbf x}$ and ${\mathbf y}$ to form the respective analogues
of (\ref{a3*}) and (\ref{a3}) and denote them by $\bbS$ and $\bbT$,
respectively. Suppose that the underlying random variables involved in
$\bbS_{{\mathbf x}{\mathbf y}},\bbT_{{\mathbf x}{\mathbf y}}$ and $\bbS,\bbT$ are
in the same probability space $(\Omega,P)$.

Recall the definitions of $G^{(i)}_{p_1,p_2},i=1,2$ in (\ref{c1}) and
(\ref{c2}), and let
$
R_n^{(i)}=\int\lambda \,dG^{(i)}_{p_1,p_2}.
%&&R_n^{(2)}=p_1\int\lambda \,d\Big(F^{\bbT_{xy}}(
%D^{(i)}=p_1\int\lambda \,d\Big(F^{\bbR_{xy}^{(i)}}_{\mathbb{H}_1}(
$
%where $\bbR_{xy}^{(1)}$ represents the matrix $\bbS_{\bbx\bby}$ while $
%$F^{\bbR_{xy}^{(i)}}_{\mathbb{H}_0},F^{\bbR_{xy}^{(i)}}_{
%while $F^{\bbR_{xy}^{(i)}}_{\mathbb{H}_1}$ is the ESD of
%$\bbR_{xy}^{(i)}$ under $\mathbb{H}_1$.

%th4.1 #&#
\begin{teo}\label{teo4}
In addition to assumptions in Theorems~\ref{teo1}~or~\ref{teo3} suppose that for any $M>0$
%
%e4.1 #&#
\begin{equation}
\label{w4} P \bigl(\bigl\llvert \tr (\bbS-\bbS_{{\mathbf x}{\mathbf y}})\bigr
\rrvert
\geq M \bigr)\rightarrow1,\qquad P \bigl(\bigl\llvert \tr (\bbT-\bbT_{{\mathbf
x}{\mathbf y}})
\bigr\rrvert\geq M \bigr)\rightarrow1.
\end{equation}
%
%where $\bbS_{\bbx\bby}^{\mathbb{H}_j}$ is $\bbS_{\bbx\bby}$ under $
%defined similarly.
Then
%
%e4.2 #&#
\begin{equation}
\lim_{n\rightarrow\infty}P\bigl(R^{(i)}_n>z^{(i)}_{1-\alpha}
\mbox{ or } R^{(i)}_n<z^{(i)}_\alpha\mid
\mathbb{H}_1\bigr)=1,
\end{equation}
where $z^{(i)}_{1-\alpha}$ and $z^{(i)}_{\alpha}$ are, respectively,
$(1-\alpha)$ and $\alpha$ quantiles of the asymptotic distribution of
the statistic $R_n^{(i)}$ under the null hypothesis.
\end{teo}

%re7 #&#
\begin{rmk}
For example, one may take $\bbS=(\bbX\bbL\bbX^T)^{-1}\bbX\bbL\bbY
^T\times\break (\bbY\bbL\bbY^T)^{-1}\bbY\bbL\bbX^T$ and $\bbS_{{\mathbf x}{\mathbf
y}}=(\bbX\bbX^T)^{-1}\bbX\bbP_{{\mathbf y}}\bbX^T$ with $\bbL$ being
a random matrix and $\bbP_{{\mathbf y}}=\bbY^{T}(\bbY\bbY^{T})^{-1}\bbY
$. Particularly, if $\bbL=\bbI+{\mathbf e}{\mathbf e}^T$ with ${\mathbf
e}=x^2(1,1,\ldots,1)$ and $x^2$ having finite moment,
then under assumptions in Theorems~~\ref{teo1}~or~\ref{teo3} it
can be proved that
\[
\tr  (\bbS-\bbS_{{\mathbf x}{\mathbf y}} )=O_p(n)
\]
satisfying (\ref{w4}).
\end{rmk}

Next, we evaluate the powers of the modified statistics with the
dividing-sample method.
Draw $n$ samples from alternatives ${\mathbf x}$ and ${\mathbf y}$ to form the
respective analogues of $\bbT_{{\mathbf x}{\mathbf y}}^{(i)}$, $i=1,2$,
$\widetilde{\bbT}_{{\mathbf x}{\mathbf y}}^{(1)}$ and denote them by $\bbT
^{(i)}$, $i=1,2$, $\widetilde{\bbT}{}^{(1)}$, respectively.
Let
\begin{eqnarray*}
J_n^{(1)}&=&\int\lambda \,dF^{\bbT^{(1)}}(\lambda)-2\int
\lambda \,dF^{\bbT^{(2)}}(\lambda),
\nonumber
\\
J_n^{(2)}&=&\int\lambda \,dF^{\bbT^{(1)}}(\lambda)-\int
\lambda \,dF^{\bbT^{(2)}}(\lambda),
\nonumber
\\
J_n^{(3)}&=&\int\lambda \,dF^{\widetilde{\bbT}{}^{(1)}}(\lambda)-\int
\lambda \,dF^{\bbT^{(2)}}(\lambda).
\end{eqnarray*}

%th4.2 #&#
\begin{teo}\label{teo7}
In addition to assumptions in Theorem~\ref{teo3}, suppose that for any $M>0$,
%
%e4.3 #&#
%e4.4 #&#
%e4.5 #&#
%e4.6 #&#
%e4.7 #&#
%e4.8 #&#
\begin{eqnarray}
P \bigl(\bigl\llvert \tr \bigl(\bbT^{(1)}\bigr)-2\tr \bigl(
\bbT^{(2)}\bigr)- \bigl(\tr \bigl(\bbT_{{\mathbf
x}{\mathbf y}}^{(1)}
\bigr)-2\tr \bigl(\bbT_{{\mathbf x}{\mathbf y}}^{(2)}\bigr) \bigr)\bigr\rrvert
\geq
M \bigr)&\rightarrow&1,
\nonumber\\[-8pt] \label{add60y}\\[-8pt]
\eqntext{\mbox{if }c_2'<1/2;}
\\
P \bigl(\bigl\llvert \tr \bigl(\bbT^{(1)}\bigr)-\tr \bigl(
\bbT^{(2)}\bigr)- \bigl(\tr \bigl(\bbT_{{\mathbf
x}{\mathbf y}}^{(1)}
\bigr)-\tr \bigl(\bbT_{{\mathbf x}{\mathbf y}}^{(2)}\bigr) \bigr)\bigr\rrvert
\geq M
\bigr)&\rightarrow&1,
\nonumber\\[-8pt] \label{add61y}\\[-8pt]
\eqntext{\mbox{if }c_2'\geq1;}
\\
P \bigl(\bigl\llvert \tr \bigl(\widetilde{\bbT}{}^{(1)}\bigr)-\tr \bigl(
\bbT^{(2)}\bigr)- \bigl(\tr \bigl(\widetilde{\bbT}_{{\mathbf x}{\mathbf y}}^{(1)}
\bigr)-\tr \bigl(\bbT_{{\mathbf x}{\mathbf
y}}^{(2)}\bigr) \bigr)\bigr\rrvert\geq M
\bigr)&\rightarrow&1,
\nonumber\\[-8pt] \label{add62y} \\[-8pt]
\eqntext{\mbox{if }1/2\leq c_2'<1.}
\end{eqnarray}
Then
%if $c_2'<1/2$,
%
\[
\lim_{n\rightarrow\infty}P\bigl(J_n^{(i)}>z_{1-\alpha}^{(i)}
\mbox{ or }J_n^{(i)}<z_{\alpha}^{(i)}\mid
\mathbb{H}_1\bigr)=1,\qquad i=1,2,3,
\]
%
%if $c_2'\geq1$,
%%J_n^{(2)}<z_{\alpha}^{(2)}\right\vert\mathbb{H}_1)=1;
%if $1/2\leq c_2'<1$,
%%J_n^{(3)}<z_{\alpha}^{(3)}\left\vert\mathbb{H}_1)=1,
where $z^{(i)}_{1-\alpha}$ and $z^{(i)}_{\alpha}$ are, respectively,
$(1-\alpha)$ and $\alpha$ quantiles of the asymptotic distribution of
the statistic $J_n^{(i)}$ under the null hypothesis, $i=1,2,3$.
\end{teo}

%s5 #&#
\section{Applications of CCA}\label{sec5}
%Canonical correlation analysis has been widely used in the literature
%to investigate the relationship between two random variable sets.
%For any two high dimensional random vectors $\bbx\dvtx  p_1\times1$ and $
%and \bby are dependent.
%Based on the CLT developed in the last section, we propose the
%following statistic
%S_n=\sum^{p_1}_{i=1}r_i^{2},
%where $r_i^{2}, i=1,2,\ldots,p_1$ are the eigenvalues of the matrix $

%Note that
%S_n=\int\lambda \,dG_{p_1,p_2}(\lambda).
%Hence the CLT developed in Theorem \ref{teo1} provides the null
%asymptotic distribution for the independence test (

%Independence tests for two different random variable sets play an
%important role in statistical inference. It is often the statistical
%tool used to identify the simplifying structure of the data. The
%following two examples are two important problems in multivariate
%analysis and time series respectively. They can be transformed into
%independence test based on canonical correlation analysis and then we
%can utilize the proposed test statistic to solve these problem for
%high dimensional data.
This section explores some applications of the proposed test. We
consider two examples from multivariate analysis and time series
analysis, respectively.

%s5.1 #&#
\subsection{Multivariate regression test with CCA}\label{sec5.1}

Consider the multivariate regression (MR) model as follows:
%
%e5.1 #&#
\begin{equation}
\label{model1} \bbY=\bbX\bbB+\bbE,
\end{equation}
where
\begin{eqnarray*}
\bbY&=&[{\mathbf y}_1,{\mathbf y}_2,\ldots,{\mathbf
y}_{p_1}]_{n\times p_1},\qquad \bbX=[\bolds{1}_n,{\mathbf
x}_1,{\mathbf x}_2,\ldots,{\mathbf x}_{p_2}]_{n\times p_2},
\\
\bbB&=&[\bolds{\beta}_1, \bolds{\beta}_2,
\ldots, \bolds{\beta}_{p_1}]_{p_2\times p_1},\qquad \bbE=[{\mathbf
e}_1,{\mathbf e}_2,\ldots,{\mathbf e}_{p_1}]_{n\times p_1},
\end{eqnarray*}
and each of the vectors ${\mathbf y}_j$, ${\mathbf x}_j$, ${\mathbf e}_j$, for
$j=1,2,\ldots,p_1$ is $n\times1$ vectors and $\{\bolds{\beta
}_i$, $i=1,2,\ldots,p_1\}$ are $p_2\times1$ vectors.

Let $\bbA_{{\mathbf x}{\mathbf y}}=\frac{1}{n}\bbX^{T}\bbY$ and $\bbA
_{{\mathbf x}{\mathbf x}}=\frac{1}{n}\bbX^{T}\bbX$. We have the least square
estimate of $\bbB$
%
%e5.2 #&#
\begin{equation}
\hat\bbB=\bbA_{{\mathbf x}{\mathbf x}}^{-1}\bbA_{{\mathbf x}{\mathbf y}}.
\end{equation}
The most common hypothesis testing is to test whether there exists
linear relationship between the two sets of variables (response
variables and predictor variables) or the overall regression test
%
%e5.3 #&#
\begin{equation}
\mathbb{H}_0\dvtx  \bbB=\bolds{0}.
\end{equation}
%
%Selecting $\bbC_{k\times q}=[\bolds{0}, \bbI_k]$ of full rank $k$ and
%$\bbM_{p\times p}=\bbI_p$, the test that $\bbB_1=\bolds{0}$ is easily
%derived from the general matrix form of the hypothesis, $\bbC\bbB\bbM=

To test $\mathbb{H}_0\dvtx  \bbB=\bolds{0}$, Wilks' $\Lambda$ criterion is
%
%e5.4 #&#
\begin{equation}
\label{a2} \Lambda=\frac{\det(\bbE)}{\det(\bbE+\bbH)}=\prod^{s}_{i=1}(1+
\lambda_i)^{-1},
\end{equation}
where
%
%e5.5 #&#
\begin{equation}
\bbE=\bbY^{T} \bigl(\bbI-\bbX\bigl(\bbX^{T}\bbX
\bigr)^{-1}\bbX^{T} \bigr)\bbY
\end{equation}
and
%
%e5.6 #&#
\begin{equation}
\bbH=\hat\bbB^{T}\bigl(\bbX^{T}\bbX\bigr)\hat\bbB;
\end{equation}
and $\{\lambda_i\dvtx  i=1,\ldots,s\}$ are the roots of $\det(\bbH-\lambda
\bbE)=0$, $s=\min(k,p)$. An alternative form for $\Lambda$ is to
employ sample covariance matrices. That is, $\bbH=\bbA_{{\mathbf y}{\mathbf
x}}\bbA_{{\mathbf x}{\mathbf x}}^{-1}\bbA_{{\mathbf x}{\mathbf y}}$ and $\bbE=\bbA
_{{\mathbf y}{\mathbf y}}-\bbA_{{\mathbf y}{\mathbf x}}\bbA_{{\mathbf x}{\mathbf x}}^{-1}\bbA
_{{\mathbf x}{\mathbf y}}$, so that $\det(\bbH-\lambda\bbE)=0$ becomes
$\det (\bbA_{{\mathbf y}{\mathbf x}}\bbA_{{\mathbf x}{\mathbf x}}^{-1}\bbA_{{\mathbf
x}{\mathbf y}}-\lambda(\bbA_{{\mathbf y}{\mathbf y}}-\bbA_{{\mathbf y}{\mathbf x}}\bbA
_{{\mathbf x}{\mathbf x}}^{-1}\bbA_{{\mathbf x}{\mathbf y}}) )=0$. From Theorem~2.6.8 of \cite{NT2001}, we have $\det(\bbH-\theta(\bbH+\bbE
))=\det (\bbA_{{\mathbf y}{\mathbf x}}\bbA_{{\mathbf x}{\mathbf x}}^{-1}\bbA
_{{\mathbf x}{\mathbf y}}-\theta\bbA_{{\mathbf y}{\mathbf y}} )=0$ so that
%
%e5.7 #&#
\begin{equation}
\Lambda=\prod^{s}_{i=1}(1+
\lambda_i)^{-1}=\prod^{s}_{i=1}(1-
\theta_i)=\frac{\det(\bbA_{{\mathbf y}{\mathbf y}}-\bbA_{{\mathbf y}{\mathbf x}}\bbA
_{{\mathbf x}{\mathbf x}}^{-1}\bbA_{{\mathbf x}{\mathbf y}})}{\det(\bbA_{{\mathbf y}{\mathbf y}})}.
\end{equation}
%
%From this, we can see that testing $\mathbb{H}_0\dvtx  \bbB_1=\bolds{0}$
%is equivalent to testing $\Sigma_{\bbx\bby}=\bolds{0}$.
Evidently, the quantities $r_i^{2}=\theta_i, i=1,\ldots,s$ are sample
canonical correlation coefficients. Therefore, the test statistic (\ref
{a2}) can be rewritten as
%
%e5.8 #&#
\begin{equation}
\log\Lambda=\sum^{s}_{i=1}\log
\bigl(1-r_i^2\bigr).
\end{equation}
From this point of view, the multiple regression test is equivalent to
the independence test based on canonical correlation coefficients.
As stated in the last section, the statistic $\log\Lambda$ is not
stable in the high dimensional cases. Hence, our test statistic $S_n$
or $T_n$ can be applied in the MR test.

%s5.2 #&#
\subsection{Testing for cointegration with CCA}\label{sec5.2}

Consider an $n$-dimensional vector process $\{{\mathbf y}_t\}$ that has a
first-order error correction representation
%
%e5.9 #&#
\begin{equation}
\Delta{\mathbf y}_t=-\bolds{\alpha}\bolds{
\beta}'{\mathbf y}_{t-1}+\bolds{\varepsilon}_t,\qquad
t=1,\ldots,T,
\end{equation}
where $\bolds{\alpha}$ and $\bolds{\beta}$ are full rank
$n\times r$ matrices $(r<n)$ and the $n$-dimensional innovation $\{
\varepsilon_t\}$ is i.i.d. with zero mean and positive covariance
matrix $\bolds{\Omega}$. Select $\bolds{\alpha}$ and
$\bolds{\beta}$ so that the fact that $\llvert \bbI_n-(\bbI
_n-\bolds{\alpha}\bolds{\beta}')z\rrvert =0$ implies that
either $\llvert z\rrvert >1$ or $z=1$ and that $\bolds{\alpha
}_{\bot
}'\bolds{\beta}_{\bot}$ is of full rank, where $\bolds
{\alpha}_{\bot}$ and $\bolds{\beta}_{\bot}$ are full rank
$n\times(n-r)$ matrices orthogonal to $\bolds{\alpha}$ and
$\bolds{\beta}$. Under these assumptions, $\{{\mathbf y}_t\}$ is
$I(1)$ with $r$ cointegration relations among its elements; that is, $\{
\bolds{\beta}'{\mathbf y}_t\}$ is $I(0)$. Here, $I(d)$ denotes
integrated of order $d$.

The goal is to test
%
%e5.10 #&#
\begin{equation}
\label{cointegrationtest} \mathbb{H}_0\dvtx  r=0~(\bolds{\alpha
}=\bolds{
\beta}=\bolds{0});\quad\mbox{against}\quad\bbH_1\dvtx  r>0;
\end{equation}
that is, whether there exists cointegration relationships among the
elements of the time series $\{{\mathbf y}_t\}$.\vadjust{\goodbreak}

This cointegration test is equivalent to testing:
\begin{eqnarray*}
&&\mathbb{H}_0\dvtx  \Delta{\mathbf y}_t\mbox{ is independent with }
\Delta{\mathbf y}_{t-1};\quad\mbox{against}
\\
&& H_1\dvtx  \Delta{\mathbf
y}_t\mbox{ is dependent with }\Delta{\mathbf y}_{t-1}.
\end{eqnarray*}

In order to apply canonical correlation coefficients to cointegration
test (\ref{cointegrationtest}), we construct random matrices
%
%e5.11 #&#
%e5.12 #&#
\begin{eqnarray}
\label{cox} \bbX&=& (\Delta{\mathbf y}_2, \Delta{\mathbf y}_4,
\ldots, \Delta{\mathbf y}_{2t-2}, \Delta{\mathbf y}_{2t}, \ldots,
\Delta{\mathbf y}_{T} ),
\\
\label{coy} \bbY&=& (\Delta{\mathbf y}_{1}, \Delta{\mathbf
y}_{3}, \ldots, \Delta{\mathbf y}_{2t-1}, \Delta{\mathbf
y}_{2t+1}, \ldots, \Delta{\mathbf y}_{T-1} ).
\end{eqnarray}

%s6 #&#
\section{Simulation results}\label{sec6}

This section reports some simulated examples to show the finite sample
performance of the
proposed test.

%s6.2 #&#
\subsection{Empirical sizes and empirical powers}\label{sec6.2}

First, we introduce the method of calculating empirical sizes and
empirical powers. Let $z_{1-\alpha}$ be the $100(1-\alpha)\%$
quantile of the asymptotic null distribution of the test statistic
$S_n$. With $K$ replications of the data set simulated under the null
hypothesis, we calculate the empirical size as
%
%e6.2 #&#
\begin{equation}
\hat\alpha=\frac{\{\sharp\mbox{ of }S_n^H\geq z_{1-\alpha}\}}{K},
\end{equation}
where $S_n^H$ represents the values of the test statistic $S_n$ based
on the data simulated under the null hypothesis.

The empirical power is calculated as
%
%e6.3 #&#
\begin{equation}
\hat\beta=\frac{\{\sharp\mbox{ of }S_n^A\geq\hat z_{1-\alpha}\}}{K},
\end{equation}
where $S_n^A$ represents the values of the test statistic $S_n$ based
on the data simulated under the alternative hypothesis.

In our simulations, we choose $K=1000$ as the number of repeated
simulations. The significance level is $\alpha=0.05$.

%s6.3 #&#
\subsection{Testing independence}\label{sec6.3}
Consider the data generating process
%
%e6.4 #&#
\begin{equation}
{\mathbf x}=\bolds{\Sigma}_{{\mathbf x}{\mathbf x}}^{1/2}{\mathbf w},\qquad {\mathbf y}=
\bolds{\Sigma}_{{\mathbf y}{\mathbf y}}^{1/2}{\mathbf v},
\end{equation}
%
%and two cases are investigated as
with
\begin{eqnarray*}
\mbox{(a)}\quad \bolds{\Sigma}_{{\mathbf x}{\mathbf x}}&=&\bbI_{p_1},\qquad \bolds{
\Sigma}_{{\mathbf y}{\mathbf y}}=\bbI_{p_2};
\\
\mbox{(b)}\quad \bolds{\Sigma}_{{\mathbf x}{\mathbf x}}&=&\bigl(\sigma_{kh}^{SP}
\bigr)^{p_1}_{k,h=1},\qquad \bolds{\Sigma}_{{\mathbf y}{\mathbf
y}}=\bbI_{p_2};
\\
\mbox{(c)}\quad \bolds{\Sigma}_{{\mathbf x}{\mathbf x}}&=&\bigl(\sigma_{kh}^{AR}
\bigr)^{p_1}_{k,h=1},\qquad \bolds{\Sigma}_{{\mathbf y}{\mathbf
y}}= \bbI_{p_2};
\\
\mbox{(d)}\quad \bolds{\Sigma}_{{\mathbf x}{\mathbf
x}}&=&\bbB'\operatorname{cov}(\bbf_t)\bbB+\bolds{\Sigma}_{\bbu},\qquad
\bolds{\Sigma}_{{\mathbf y}{\mathbf y}}=\bbI_{p_2},
\end{eqnarray*}
where
\begin{eqnarray*}
\sigma^{SP}_{kh}=\cases{ 1, &\quad$k=h$;
\cr
\theta, &\quad
$k=1$; $h=2,3,\ldots,\bigl[p_1^{1/3}\bigr]$; or $h=1$;
$k=2,3,\ldots,\bigl[p_1^{1/3}\bigr]$;
\cr
0, &\quad others}
\end{eqnarray*}
with $\theta=0.2$ and
\[
\sigma_{kh}^{AR}=\frac{\phi^{\llvert k-h\rrvert }}{1-\phi^2},\qquad
k,h=1,2,\ldots,p_1, \phi=0.8.
\]
Here, $\bbB=\frac{1}{\sqrt{p_1}}({\mathbf b}_1,{\mathbf b}_2,\ldots,{\mathbf
b}_{p_1})$ is a deterministic matrix. In the simulation, each ${\mathbf
b}_i\dvtx  r\times1$\vadjust{\goodbreak} is generated independently from a normal distribution
with covariance matrix being an $r\times r$ identity matrix and mean
$\bolds{\mu}_B$ consisting of all $1$.
$\operatorname{cov}(\bbf_t)$ is also an $r\times r$ identity matrix and $\bolds
{\Sigma}_{\bbu}$ is a $p_1\times p_1$ identity matrix.

%t1 #&#
\begin{table}
\tabcolsep=0pt
\caption{Empirical sizes of the proposed test $S_n$ and the renormalized
likelihood ratio test $\mathrm{MLR}_n$ at $0.05$ significance level for DGP~\textup{(a)} and DGP~\textup{(b)}}\label{tb1}
\begin{tabular*}{\tablewidth}{@{\extracolsep{\fill}}@{}lcccc@{}}
  \hline
  $\bolds{(p_1,p_2,n)}$ & $\bolds{S_n}$ \textbf{DGP~(a)}  & $\bolds{S_n}$ \textbf{DGP~(b)} & $\mathbf{MLR}_{\bolds{n}}$ \textbf{DGP~(a)}
  & $\mathbf{MLR}_{\bolds{n}}$ \textbf{DGP~(b)} \\
\hline
  $(10,20,40)$ & 0.0458 & 0.0461 & 0.0481 & 0.0490\\
  $(20,30,60)$ & 0.0480 & 0.0488 & 0.0440 & 0.0448\\
  $(30,60,120)$ & 0.0475 & 0.0480 & 0.0530 & 0.0520 \\
  $(40,80,160)$ & 0.0464 & 0.0466 & 0.0420 & 0.0420\\
  $(50,100,200)$ & 0.0503 & 0.0504 & 0.0487 & 0.0500  \\
  $(60,120,240)$ & 0.0490 & 0.0490 & 0.0574 & 0.0572\\
  $(70,140,280)$ & 0.0524 & 0.0520 & 0.0570 & 0.0582\\
  $(80,160,320)$ & 0.0500 & 0.0500 & 0.0632 & 0.0583\\
  $(90,180,360)$ &  0.0521 & 0.0511 & 0.0559 & 0.0580\\
  $(100,200,400)$ & 0.0501 & 0.0503 & 0.0482 & 0.0589 \\
  $(110,220,440)$ & 0.0504 & 0.0500 & 0.0440 & 0.0590\\
  $(120,240,480)$ & 0.0513 & 0.0511 & 0.0400 & 0.0432\\
  $(130,260,520)$ & 0.0511 & 0.0511 & 0.0520 & 0.0560\\
  $(140,280,560)$ & 0.0469 & 0.0474 & 0.0582 & 0.0580\\
  $(150,300,600)$ & 0.0495 & 0.0500 & 0.0590 & 0.0593\\
  $(160,320,640)$ & 0.0514 & 0.0517 & 0.0437 & 0.0559\\
  $(170,340,680)$ & 0.0498  & 0.0500 & 0.0428 & 0.0430\\
  $(180,360,720)$ & 0.0509 & 0.0510 & 0.0580 & 0.0577\\
  $(190,380,760)$ & 0.0488 & 0.0485 & 0.0388 & 0.0499 \\
  $(200,400,800)$ & 0.0491 & 0.0491 & 0.0462 & 0.0499\\
  $(210,420,840)$ & 0.0491  & 0.0500 & 0.0450 & 0.0555\\
  $(220,440,880)$ & 0.0515 & 0.0510 & 0.0572 & 0.0588\\
  $(230,460,920)$ & 0.0493 & 0.0498 & 0.0470 & 0.0488\\
  $(240,480,960)$ & 0.0482 & 0.0479 & 0.0521 & 0.0561\\
  $(250,500,1000)$ & 0.0452 & 0.0450 & 0.0527 & 0.0545\\
  \hline
\end{tabular*}\vspace*{-6pt}
\end{table}

%t2 #&#
\begin{table}%[h]
\tabcolsep=0pt
\caption{Empirical sizes of the proposed test $T_n$ at $0.05$ significance level for DGP~\textup{(a)}--DGP~\textup{(d)}}\label{tb11}
\begin{tabular*}{\tablewidth}{@{\extracolsep{\fill}}@{}lcccc@{}}
  \hline
  $\bolds{(p_1,p_2,n)}$ & $\bolds{T_n}$ \textbf{DGP~(a)}  & $\bolds{T_n}$ \textbf{DGP~(b)}
  & $\bolds{T_n}$ \textbf{DGP~(c)}  & $\bolds{T_n}$ \textbf{DGP~(d)}\\
\hline
  $(100,50,80)$ & 0.0569 & 0.0462 & 0.0642  & 0.0602\\
  $(140,70,120)$ & 0.0573 & 0.0429 & 0.0619 & 0.0600\\
  $(180,90,150)$ &  0.0577 & 0.0452 & 0.0623  & 0.0583 \\
  $(200,100,170)$ & 0.0552 & 0.0429 & 0.0594 & 0.0592  \\
  $(240,120,180)$ & 0.0581 & 0.0510 & 0.0602 & 0.0608 \\
  $(280,140,250)$ & 0.0571 & 0.0483 & 0.0592 & 0.0584\\
  $(320,160,270)$ & 0.0521 & 0.0479 & 0.0603 & 0.0549\\
  $(360,180,290)$ & 0.0529 & 0.0489 & 0.0574 & 0.0569\\
  $(400,190,300)$ & 0.0542 & 0.0522 & 0.0589 & 0.0579 \\
  $(440,220,330)$ & 0.0557 & 0.0529 & 0.0542  & 0.0581 \\
  $(480,240,350)$ & 0.0531 & 0.0562 & 0.0579 & 0.0569\\
  \hline
\end{tabular*}
\tabnotetext[]{ta2}{The parameter $t$ in the statistic $T_n$ takes a value of 40.
For DGP~(a), we use the original statistic $T_n$ in Theorem \ref{teo3};
for DGP~(b), the statistic in Theorem \ref{teo5} is used; for DGPs (c) and (d), the dividing-sample statistic in Theorem \ref{teo6} is utilized.}
\end{table}

The empirical sizes of the proposed statistics $S_n$ for data
generating processes (DGPs) (a) and (b) are listed in Table~\ref{tb1}.
Moreover, the empirical sizes for the renormalized statistic $\MLR_n$
are included as comparison with $S_n$. Here the renormalized statistic
$\MLR_n$ means the statistic
\[
p_1\int\log(1-\lambda)\,d \bigl(F^{\bbS_{{\mathbf x}{\mathbf y}}}(\lambda
)-F^{c_{1n},c_{2n}}(\lambda) \bigr).
\]
The empirical sizes of $T_n$ for DGPs (a)--(d) are listed in Table~\ref
{tb11}. For DGP~(a), we use the original statistic $T_n$; for DGP~(b),
the statistic in Theorem \ref{teo5} is used; for DGPs (c) and (d),
the dividing-sample statistic in Theorem \ref{teo6} is utilized. For
Theorem \ref{teo6}, we follow the implementation details in Appendix
of \cite{Karoui2008} to estimate the LSD $H$ for the matrix
$\bolds{\Sigma}_{{\mathbf x}{\mathbf x}}$.\vadjust{\goodbreak}

From the results in Tables~\ref{tb1} and \ref{tb11}, the proposed
statistics $S_n$ and $T_n$ work well under Assumptions \ref{ass1} and \ref{ass2}, respectively.

%%The covariance matrix $\bolds{\Sigma}_{\bbx\bbx}$ in case (b) is
%and we use the thresholding method in \cite{BL2008} to estimate it.
%A banded type matrix in (c) and a sparse matrix in (b) are both
%estimated by the thresholding method in \cite{BL20081}. A low rank
%matrix plus a sparse matrix in (d) is estimated by combining principle
%component analysis and thresholding method originated in
%
%%Moreover, for the proposed test statistic $T_n$, we use
%dividing-sample method to do the test for cases (c) and (d).

%s6.4 #&#
\subsection{Factor model dependence}\label{sec6.4}\label{factor2}
We consider the factor model as follows:
%
%e6.5 #&#
\begin{eqnarray}
\label{factor1} {\mathbf x}_t=\Lambda_1\bbf_t+
\bbu_t,\qquad {\mathbf y}_t=\Lambda_2\bbf
_t+{\mathbf v}_t, t=1,2,\ldots,n,
\end{eqnarray}
where $\Lambda_1=\frac{1}{\sqrt{p_1}}(\bolds{\lambda
}_1^{(1)}, \bolds{\lambda}_2^{(1)},\ldots, \bolds{\lambda
}_{r}^{(1)})$ and $\Lambda_2=\frac{1}{\sqrt{p_2}}(\bolds
{\lambda}_1^{(2)}, \bolds{\lambda}_2^{(2)}, \ldots,
\bolds{\lambda}_{r}^{(2)})$ are\break \mbox{$p_1\times r$} and $p_2\times r$
deterministic matrices, respectively. In the simulation, all the
components of $\bolds{\lambda}_k^{(j)}\dvtx  k=1,2,\ldots,r; j=1,2$
are generated from a normal distribution with mean being $0.8$ and
variance being $1$. $\bbf_t, t=1,2,\ldots,n$ are $r\times1$ random
vectors with i.i.d. standard Gaussian distributed elements and $\bbu
_t$ and ${\mathbf v}_t, t=1,2,\ldots,n$ are independent random vectors
whose elements are all standard Gaussian distributed.

%t3 #&#
\begin{table}[t]
\tabcolsep=0pt
\caption{Empirical powers of the proposed test $S_n$ at $0.05$ significance level for factor models}\label{tb2}
\begin{tabular*}{\tablewidth}{@{\extracolsep{\fill}}@{}lcccc@{}}
  \hline
  $\bolds{(p_1,p_2,n)}$ & $\bolds{r=3}$ & $\bolds{r=5}$ & $\bolds{r=7}$ & $\bolds{r=10}$ \\
  \hline
  $(10,20,40)$ & 0.3750 & 0.5200 & 0.5910 & 0.9320\\
  $(30,60,120)$ & 0.3070 & 0.6240 & 0.8450 & 0.9330\\
  $(50,100,200)$ & 0.3090 & 0.6700 & 0.7980 & 0.9700\\
  $(70,140,280)$ & 0.3520 & 0.6470 & 0.8330 & 0.9850\\
  $(90,180,360)$ &  0.3670 & 0.6720 & 0.8230 & 0.9880\\
  $(110,220,440)$ & 0.3570 & 0.6690 & 0.8490 & 0.9850\\
  $(130,260,520)$ & 0.3440 & 0.6390 & 0.8510 & 0.9960\\
  $(150,300,600)$ & 0.3780 & 0.6440 & 0.8370 & 0.9990\\
  $(170,340,680)$ & 0.3580  & 0.6580 & 0.8590 & 1.0000\\
  $(190,380,760)$ & 0.3490 & 0.6620 & 0.8720 & 1.0000\\
  $(210,420,840)$ & 0.3460  & 0.6790 & 0.8890 & 1.0000\\
  $(230,460,920)$ & 0.3800 & 0.6930 & 0.8770 & 1.0000\\
  $(250,500,1000)$ & 0.3470 & 0.6890 & 0.8940 & 1.0000\\
  \hline
\end{tabular*}
\tabnotetext[]{ta3}{The powers are under the alternative hypothesis that $\bbx$ and $\bby$ satisfy the factor model (\ref{factor1}). $r$ is the number of factors.}\vspace*{-7pt}
\end{table}

%t4 #&#
\begin{table}[b]
\tabcolsep=0pt
\caption{Empirical powers of the proposed test $T_n$ at $0.05$ significance level for factor models}\label{tb22}
\begin{tabular*}{\tablewidth}{@{\extracolsep{\fill}}lcccc@{}}
  \hline
  $\bolds{(p_1,p_2,n)}$ & $\bolds{r=3}$ & $\bolds{r=5}$ & $\bolds{r=7}$ & $\bolds{r=10}$ \\
  \hline
  $(100,50,80)$ & 0.3680 & 0.6380 & 0.7330 & 0.9470\\
  $(140,70,120)$ & 0.3380 & 0.6440 & 0.8690 & 0.9520\\
  $(180,90,150)$ & 0.3290 & 0.6190 & 0.8890 & 0.9740\\
  $(200,100,170)$ & 0.3410 & 0.6270 & 0.8920 & 0.9820\\
  $(240,120,180)$ &  0.3340 & 0.6290 & 0.8840 & 0.9790\\
  $(280,140,250)$ & 0.3570 & 0.6480 & 0.8730 & 0.9870\\
  $(320,160,270)$ & 0.3490 & 0.7120 & 0.8890 & 0.9940\\
  $(360,180,290)$ & 0.3690 & 0.6890 & 0.8930 & 0.9920\\
  $(400,200,310)$ & 0.3830  & 0.7080 & 0.9030 & 0.9980\\
  $(440,220,330)$ & 0.3920 & 0.7040 & 0.8930 & 1.0000\\
  $(480,240,350)$  & 0.3970 & 0.6990 & 0.9110 & 1.0000\\
  \hline
\end{tabular*}
\tabnotetext[]{ta4}{The powers are under the alternative hypothesis that $\bbx$ and $\bby$
satisfy the factor model (\ref{factor1}). $r$ is the number of factors. The parameter $t$ in the statistic $T_n$ takes value of $40$. For $T_n$, we use its modified dividing-sample version in Theorem \ref{teo6}.}
\end{table}

%t7 #&#
\begin{table}%[h]
\tabcolsep=0pt
\tablewidth=220pt
\caption{Empirical powers of the proposed test $S_n$ at $0.05$ significance level for uncorrelated but dependent case}\label{tb4}
\begin{tabular*}{\tablewidth}{@{\extracolsep{\fill}}@{}lcc@{}}
  \hline
  $\bolds{(p_1,p_2,n)}$ & $\bolds{\omega=4}$  & $\bolds{\omega=10}$ \\
\hline
  $(10,20,40)$ & 0.8140 & 0.9690 \\
  $(30,60,120)$ & 0.8200 & 0.9510 \\
  $(50,100,200)$ & 0.8220 & 0.9600 \\
  $(70,140,280)$ & 0.8100 & 0.9610 \\
  $(90,180,360)$ &  0.8210 & 0.9640 \\
  $(110,220,440)$ & 0.8110 & 0.9670 \\
  $(130,260,520)$ & 0.8320 & 0.9740 \\
  $(150,300,600)$ & 0.8420 & 0.9740 \\
  $(170,340,680)$ & 0.8450  & 0.9760 \\
  $(190,380,760)$ & 0.8580 & 0.9680 \\
  $(210,420,840)$ & 0.8420  & 0.9670 \\
  $(230,460,920)$ & 0.8440 & 0.9810 \\
  $(250,500,1000)$ & 0.8620 & 0.9810 \\
  \hline
\end{tabular*}
\tabnotetext[]{ta7}{The powers are under the alternative hypothesis that $Y_{it}=X_{it}^{\omega}-EX_{it}^{\omega}, i=1,2,\ldots,p_1$ and $Y_{jt}=\varepsilon_{jt}, j=p_1+1,\ldots,p_2$; $t=1,\ldots,n$, where $\varepsilon_{jt}, j=p_1+1,\ldots,p_2$; $t=1,\ldots,n$ are standard normal distributed and independent with $X_{it}$ and $\omega=4,10$.}
\end{table}

For this model, ${\mathbf x}_t$ and ${\mathbf y}_t$ are not independent if
$r\neq0$. The proposed test statistic $S_n$ and $T_n$ can be used to
detect this dependent structure. Tables~\ref{tb2} and \ref{tb22}
illustrate the powers of the proposed statistics $S_n$ and $T_n$,
respectively, as $r$ increases from $3$ to $10$. For $T_n$, we use its
modified version in Theorem \ref{teo6}.
Results in these tables indicate that for one triple $(p_1, p_2, n)$,
the power increases as the number of factors $r$ increases.
%Meanwhile, for the fixed number of factors $r$, as $(p_1, p_2)$
%increases, the empirical powers decrease.
This phenomenon makes sense since the dependence between ${\mathbf x}_t$
and ${\mathbf y}_t$ is described by the $r$ common factors contained in the
factor vector $\bbf_t$.
%As the dimensionality $(p_1,p_2)$ tend to infinity, the dependence
%caused by the $r$ common factors becomes weak.
Stronger dependence between ${\mathbf x}_t$ and ${\mathbf y}_t$ exists while
more common factors are included in the model.

Here, we would like to point out that using CCA based on the sample
covariance matrices\vadjust{\goodbreak} with sample mean will incorrectly conclude that
${\mathbf x}_t$ and ${\mathbf y}_t$ can be independent even if $r>0$ but $\bbf
_t=\bbf$ independent of $t$ because CCA of ${\mathbf x}_t$ and ${\mathbf y}_t$
is the same as that of $\bbu_t$ and ${\mathbf v}_t$. This is why (\ref
{root1}) and (\ref{a3}) are used.

%s6.5 #&#
\subsection{Uncorrelated but dependent}\label{sec6.5}

The construction of (\ref{b3}) is based on the idea that the limit of
$F^{\bbS_{{\mathbf x}{\mathbf y}}}(x)$ could not be determined\vadjust{\goodbreak} from (\ref
{lsd}) when ${\mathbf x}$ and ${\mathbf y}$ have correlation. Thus, a natural
question is whether our statistic works in the uncorrelated but
dependent case. Below is such an example to demonstrate
the power of the test statistic in detecting uncorrelatedness.

Let ${\mathbf x}_t=(X_{1t},X_{2t},\ldots,X_{p_1t})^{T}, t=1,2,\ldots,n$
be i.i.d. normally distributed random vectors with zero means and unit
variances. Define ${\mathbf y}_t=(Y_{1t},Y_{2t},\ldots,\break  Y_{p_2t})^{T}$,
$t=1,2,\ldots,n$ by $Y_{it}=(X_{it}^{2k}-EX_{it}^{2k}), i=1,2,\ldots,
\min(p_1, p_2)$ and if $p_1<p_2$, we let $Y_{jt}=\varepsilon_{jt},
j=p_1+1,\ldots,p_2; t=1,\ldots,n$, where $\varepsilon_{jt},
j=p_1+1,\ldots,p_2; t=1,\ldots,n$ are i.i.d. normal distributed
random variables and independent with ${\mathbf x}_t$ and $k$ is an
positive integer.

%t8 #&#
\begin{table}
\tabcolsep=0pt
\tablewidth=220pt
\caption{Empirical powers of the proposed test $T_n$ at $0.05$ significance level for uncorrelated but dependent case}\label{tb44}
\begin{tabular*}{\tablewidth}{@{\extracolsep{\fill}}@{}lcc@{}}
  \hline
  $\bolds{(p_1,p_2,n)}$ & $\bolds{\omega=4}$  & $\bolds{\omega=10}$\\
\hline
  $(100,50,80)$ & 0.7240 & 0.8690\\
  $(140,70,120)$ & 0.7940 & 0.8890\\
  $(180,90,150)$ &  0.7830 & 0.8940\\
  $(200,100,170)$ & 0.7910 & 0.9340\\
  $(240,120,180)$ & 0.8420 & 0.9290\\
  $(280,140,250)$ & 0.8680 & 0.9580\\
  $(320,160,270)$ & 0.9010 & 0.9820\\
  $(360,180,290)$ & 0.9190 & 0.9940\\
  $(400,200,310)$ & 0.9530 & 0.9990\\
  $(440,220,330)$ & 0.9820 & 1.0000\\
  $(480,240,350)$ & 0.9940 & 1.0000\\
  \hline
\end{tabular*}
\tabnotetext[]{ta8}{The powers are under the alternative hypothesis that $Y_{it}=X_{it}^{\omega}-EX_{it}^{\omega}, i=1,2,\ldots,p_1$ and $Y_{jt}=\varepsilon_{jt}, j=p_1+1,\ldots,p_2$; $t=1,\ldots,n$, where $\varepsilon_{jt}, j=p_1+1,\ldots,p_2$; $t=1,\ldots,n$ are standard normal distributed and independent with $X_{it}$ and $\omega=4,10$. The parameter $t$ in the statistic $T_n$ takes value of $40$. The original statistic $T_n$ in Theorem \ref{teo3} is used.}
\end{table}

%re8 #&#
\begin{rmk}
For standard normal random variable $X_{it}$, the $2k$th moment is
$EX_{it}^{2k}=2^{-k}\frac{(2k)!}{k!}$.
\end{rmk}

For this model, ${\mathbf x}_t$ and ${\mathbf y}_t$ are uncorrelated since
$\operatorname{cov}(X_{it},Y_{it})=EX_{it}^{2k+1}-EX_{it}EX_{it}^{2k}=0$. Simulation
results in Tables~\ref{tb4}~and~\ref{tb44} provide the empirical
powers of $S_n$ and $T_n$ by taking $k=2$ and $k=5$, respectively. They
show that $S_n$ and $T_n$ can distinguish this kind of dependent
relationship well when $k=5$. For the statistic $T_n$, since the
covariance matrix of ${\mathbf x}$ is an identity matrix, we use the
original statistic $T_n$ in Theorem \ref{teo3}.

%s6.6 #&#
\subsection{ARCH type dependence}\label{sec6.6}\label{secARCH}

The statistic works in the above example because the limit of
$F^{\bbS_{{\mathbf x}{\mathbf y}}}$ cannot be determined from (\ref{lsd}) if
${\mathbf x}$ and ${\mathbf y}$ are uncorrelated. However, the limit of
$F^{\bbS_{{\mathbf x}{\mathbf y}}}(x)$ might be the same as (\ref{lsd}) when
${\mathbf x}$ and ${\mathbf y}$ are uncorrelated. We consider such an example
as follows.

Consider two random vectors ${\mathbf x}_t=(X_{1t},X_{2t},\ldots,X_{p_1t})$
and ${\mathbf y}_t=(Y_{1t},Y_{2t},\break \ldots, Y_{p_2t})$ as follows:
%
%e6.6 #&#
%e6.7 #&#
\begin{eqnarray}
\label{arch1} &\displaystyle Y_{it}=Z_{it}\sqrt{
\alpha_0+\alpha_1X_{it}^2},\qquad i=1,2,\ldots,\min(p_1, p_2);&
\\
&\displaystyle \mbox{if }p_1<p_2,\qquad Y_{jt}=Z_{jt},\qquad
j=p_1+1,\ldots,p_2,&
\end{eqnarray}
where $\bbz_t=(Z_{1t}, Z_{2t}, \ldots, Z_{p_2t})$ is a random vector
consisting of i.i.d. elements generated from Normal $(0,1)$ and $\{\bbz
_t\dvtx  t=1,\ldots,n\}$ are independent across $t$; ${\mathbf
x}_t=(X_{1t},X_{2t},\ldots,X_{p_1t})$ is also a random vector with
i.i.d. elements generated from Normal $(0,1)$ and $\{{\mathbf x}_t\dvtx  t=1,\ldots,n\}$ are independent across $t$.
Moreover, $\{\bbz_t\dvtx  t=1,\ldots,n\}$ are independent of $\{{\mathbf x}_t\dvtx  t=1,\ldots,n\}$.

For this model, ${\mathbf x}_t$ and ${\mathbf y}_t$ are dependent but
uncorrelated. Simulation results indicate that the proposed test
statistic $S_n$ cannot detect the\vspace*{1pt} dependence between them.
Nevertheless, if we substitute the elements $X_{it}^2$ and $Y_{it}^2$
for $X_{it}$ and $Y_{jt}$, respectively, in the matrix $\bbS_{{\mathbf
x}{\mathbf y}}$, then the new resulting statistic $S_n$ can capture the
dependence of this type. This efficiency is due to the correlation
between the high powers of $X_{it}$ and $Y_{it}$.

Tables~\ref{tb3} and \ref{tb33} list the powers of the proposed
statistics $S_n$ and $T_n$ for testing model (\ref{arch1}) in several
cases, that is, $\alpha_0$ and $\alpha_1$ take different values. For the
statistic $T_n$, since the covariance matrix of ${\mathbf x}$ is an
identity matrix, we use the original statistic $T_n$ in Theorem \ref
{teo3}. From the table, we can find the phenomenon that as $\alpha_1$
increases, the powers also increase. This is consistent with our
intuition because larger $\alpha_1$ brings about larger correlation
between $Y_{it}$ and $X_{it}$.

%This example makes us consider that uncorrelated relationship between $
%instead of the independence assumption. However, it is not true. The
%following example shows that uncorrelated relationship can not
%guarantee the CLT in Theorem \ref{teo1}.

%s7 #&#
\section{Empirical applications}\label{sec7}
As an application of the proposed independence test, we test the
cross-sectional dependence of daily stock returns of companies between
two different sections from New York Stock Exchange (NYSE) during the
period 2000.1.1--2002.1.1, including consumer service section,
consumer duration section, consumer nonduration section, energy
section, finance section, transport section, healthcare section,
capital goods section, basic industry section and public utility
section. The data set is obtained from Wharton Research Data Services
(WRDS) database.

We randomly choose $p_1$ and $p_2$ companies from two different
sections, respectively, such as the transport and finance section. At
each time $t$, denote the closed stock prices of these companies from
the two different sections by ${\mathbf x}_t=(x_{1t}, x_{2t}, \ldots,
x_{p_1t})^{T}$ and ${\mathbf y}_t=(y_{1t}, y_{2t}, \ldots, y_{p_2t})^{T}$,
respectively. We consider daily stock returns $\bbr_t^{{\mathbf
x}}=(r_{1t}^{{\mathbf x}}, r_{2t}^{{\mathbf x}}, \ldots, r_{p_1t}^{{\mathbf x}})$
and $\bbr_t^{{\mathbf y}}=(r_{1t}^{{\mathbf y}}, r_{2t}^{{\mathbf y}}, \ldots,
r_{p_2t}^{{\mathbf y}})$ with $r_{it}^{{\mathbf x}}=\log\frac
{x_{it}}{x_{i,t-1}}$, $i=1,2,\ldots,p_1$ and $r_{jt}^{{\mathbf y}}=\log
\frac{y_{jt}}{y_{j,t-1}}$, $j=1,2,\ldots,p_2$. The goal is to test
the dependence between $\bbr_t^{{\mathbf x}}$ and $\bbr_{t}^{{\mathbf y}}$.

%t5 #&#
\begin{table}%[h]
\tabcolsep=0pt
\caption{Empirical powers of the proposed test $S_n$ at $0.05$ significance level for $\bbx$ and $\bby$ with ARCH(1) dependent type}\label{tb3}
\begin{tabular*}{\tablewidth}{@{\extracolsep{\fill}}@{}lccccc@{}}
  \hline
  $\bolds{(p_1,p_2,n)}$ & $\bolds{(0.9,0.1)}$  & $\bolds{(0.8,0.2)}$ & $\bolds{(0.7,0.3)}$ & $\bolds{(0.6,0.4)}$ & $\bolds{(0.5,0.5)}$\\
 \hline
  $(10,20,40)$ & 0.3480 & 0.4670 & 0.6380 & 0.7650 & 0.8500\\
  $(30,60,120)$ & 0.4840 & 0.8090 & 0.9820 & 0.9990 & 1.0000\\
  $(50,100,200)$ & 0.6190 & 0.9730 & 1.0000 & 1.0000 & 1.0000\\
  $(70,140,280)$ & 0.7020 & 0.9980 & 1.0000 & 1.0000 & 1.0000\\
  $(90,180,360)$ &  0.7900 & 1.0000 & 1.0000 & 1.0000 & 1.0000\\
  $(110,220,440)$ & 0.8620 & 1.0000 & 1.0000 & 1.0000 & 1.0000\\
  $(130,260,520)$ & 0.8970 & 1.0000 & 1.0000 & 1.0000 & 1.0000\\
  $(150,300,600)$ & 0.9440 & 1.0000 & 1.0000 & 1.0000 & 1.0000\\
  $(170,340,680)$ & 0.9520  & 1.0000 & 1.0000 & 1.0000 & 1.0000\\
  $(190,380,760)$ & 0.9810 & 1.0000 & 1.0000 & 1.0000 & 1.0000\\
  $(210,420,840)$ & 0.9880  & 1.0000 & 1.0000 & 1.0000 & 1.0000\\
  $(230,460,920)$ & 0.9950 & 1.0000 & 1.0000 & 1.0000 & 1.0000\\
  $(250,500,1000)$ & 0.9980 & 1.0000 & 1.0000 & 1.0000 & 1.0000\\
  \hline
\end{tabular*}
\tabnotetext[]{ta5}{The powers are under the alternative hypothesis that $Y_{it}=Z_{it}\sqrt{\alpha_0+\alpha_1X_{it}^2}, i=1,2,\ldots,p_1; Y_{jt}=Z_{jt}, j=p_1+1,\ldots,p_2$. The pair of two numbers in this table is the value of $(\alpha_0,\alpha_1)$.}
\end{table}

%t6 #&#
\begin{table}
\caption{Empirical powers of the proposed test $T_n$ at $0.05$ significance level for $\bbx$ and $\bby$ with ARCH(1) dependent type}\label{tb33}
\begin{tabular*}{\tablewidth}{@{\extracolsep{\fill}}@{}lccccc@{}}
  \hline
  $\bolds{(p_1,p_2,n)}$ & $\bolds{(0.9,0.1)}$  & $\bolds{(0.8,0.2)}$ & $\bolds{(0.7,0.3)}$ & $\bolds{(0.6,0.4)}$ & $\bolds{(0.5,0.5)}$\\
 \hline
  $(100,50,80)$ & 0.6020 & 0.6180 & 0.7270 & 0.8930 & 0.9660\\
  $(140,70,120)$ & 0.6370 & 0.7890 & 0.8020 & 0.8990 & 0.9820\\
  $(180,90,150)$ &  0.7490 & 0.8280 & 0.9090 & 0.9920 & 1.0000\\
  $(200,100,170)$ & 0.8130 & 0.8730 & 0.9930 & 1.0000 & 1.0000  \\
  $(240,120,180)$ & 0.8920 & 0.9720 & 0.9950 & 1.0000 & 1.0000\\
  $(280,140,250)$ & 0.9470 & 0.9870 & 1.0000 & 1.0000 & 1.0000\\
  $(320,160,270)$ & 0.9900 & 0.9980 & 1.0000 & 1.0000 & 1.0000\\
  $(360,180,290)$ & 0.9910 & 0.9940 & 1.0000 & 1.0000 & 1.0000\\
  $(400,200,310)$ & 0.9890 & 0.9950 & 1.0000 & 1.0000 & 1.0000\\
  $(440,220,330)$ & 0.9920 & 1.0000 & 1.0000 & 1.0000 & 1.0000\\
  $(480,240,350)$ & 0.9980 & 0.9970 & 1.0000 & 1.0000 & 1.0000\\
  \hline
\end{tabular*}
\tabnotetext[]{ta6}{The powers are under the alternative hypothesis that $Y_{it}=Z_{it}\sqrt{\alpha_0+\alpha_1X_{it}^2}, i=1,2,\ldots,p_1; Y_{jt}=Z_{jt}, j=p_1+1,\ldots,p_2$. The pair of two numbers in this table is the value of $(\alpha_0,\alpha_1)$. The parameter $t$ in the statistic $T_n$ takes value of $40$. The original statistic $T_n$ in Theorem \ref{teo3} is used.}
\end{table}

The proposed test $S_n$ is applied to testing dependence of $\bbr
_t^{{\mathbf x}}$ and $\bbr_{t}^{{\mathbf y}}$. For each $(p_1,p_2,n)$, we
randomly choose $p_1$ and $p_2$ companies from two different sections,
construct the corresponding sample matrices $\bbX=(\bbr_1^{{\mathbf x}},
\bbr_2^{{\mathbf x}}, \ldots, \bbr_{p_1}^{{\mathbf x}})$ and $\bbY=(\bbr
_1^{{\mathbf y}}, \bbr_2^{{\mathbf y}}, \ldots, \bbr_{p_2}^{{\mathbf y}})$, and
then calculate the $P$-value by applying the proposed test. Repeat this
procedure $100$ times and derive 100 $P$-values to see whether the
cross-sectional ``dependence'' feature is popular between the tested two sections.

We test independence of daily stock returns of companies from three
pairs of sections, that is, basic industry section and capital goods
section, public utility section and capital goods section, finance
section and healthcare section. From Tables~\ref{tbaa}, \ref{tbbb}~and~\ref{tbcc}, we can see that, as the pair of numbers
of companies $(p_1,p_2)$ increases, more experiments are rejected in
terms of the $P$-values below $0.05$. It shows that cross-sectional
dependence exists and is popular for different sections in NYSE. This
suggests that the assumption that cross-sectional independence in such
empirical studies may not be appropriate.

%t9 #&#
\begin{table}
\tabcolsep=0pt
\tablewidth=242pt
\caption{$P$-values for $(p_1,p_2)$ companies from basic industry section and capital goods section of NYSE}\label{tbaa}
\begin{tabular*}{\tablewidth}{@{\extracolsep{\fill}}@{}lcc@{}}
\hline
& \multicolumn{2}{c@{}}{\textbf{No. of exp.}}\\[-6pt]
& \multicolumn{2}{c@{}}{\hrulefill}\\
\multicolumn{1}{@{}r}{$\bolds{P}$\textbf{-values:}} &  $\bolds{(p_1, p_2, n)}$ &  $\bolds{(p_1, p_2, n)}$\\
$\bolds{P}$\textbf{-value interval}  &   $\bolds{(10, 15, 20)}$ &  $\bolds{(15, 20, 25)}$\\
\hline
  $[0, 0.05]$   & 56 & 60 \\
  $[0.05, 0.1]$ & 22 & 20 \\
  $[0.1, 0.2]$  & \phantom{0}9 & 12 \\
  $[0.2, 0.3]$  & \phantom{0}2  & \phantom{0}5  \\
  $[0.3, 0.4]$  & 10  & \phantom{0}0 \\
  $[0.4, 0.5]$  & \phantom{0}1 &  \phantom{0}3\\
  $[0.6, 0.7]$  & \phantom{0}0 &  \phantom{0}0 \\
  $[0.8, 0.9]$  & \phantom{0}0 &  \phantom{0}0 \\
  $[0.9, 1]$    & \phantom{0}0 &  \phantom{0}0 \\
\hline
\end{tabular*}
\tabnotetext[]{ta9}{These are $P$-values for $(p_1,p_2)$ companies from different two sections of
NYSE: basic industry section and capital goods section,  each of which has $n$ daily stock returns during the period 2000.1.1--2002.1.1. The number of repeated experiments is $100$.
All the closed stock prices are from WRDS database. No. of Exp. is the number of experiments whose
$P$-values are in the corresponding interval.}\vspace*{-10pt}
\end{table}

%t10 #&#
\begin{table}
\tabcolsep=0pt
\tablewidth=242pt
\caption{$P$-values for $(p_1,p_2)$ companies from public utility section and capital goods section of NYSE}\label{tbbb}
\begin{tabular*}{\tablewidth}{@{\extracolsep{\fill}}@{}lcc@{}}
\hline
& \multicolumn{2}{c@{}}{\textbf{No. of exp.}}\\[-6pt]
& \multicolumn{2}{c@{}}{\hrulefill}\\
\multicolumn{1}{@{}r}{$\bolds{P}$\textbf{-values:}} &  $\bolds{(p_1, p_2, n)}$ &  $\bolds{(p_1, p_2, n)}$\\
$\bolds{P}$\textbf{-value interval}  &   $\bolds{(10, 15, 20)}$ &  $\bolds{(15, 20, 25)}$\\
\hline
  $[0, 0.05]$   & 76 & 84 \\
  $[0.05, 0.1]$ & 10  & 12 \\
  $[0.1, 0.2]$  & 4  & 2 \\
  $[0.2, 0.3]$  & 7  & 1  \\
  $[0.3, 0.4]$  & 0  & 1 \\
  $[0.4, 0.5]$  & 2  &  0\\
  $[0.6, 0.7]$  & 1  &  0 \\
  $[0.8, 0.9]$  & 0  &  0 \\
  $[0.9, 1]$    & 0  &  0 \\
  \hline
\end{tabular*}
\tabnotetext[]{ta10}{These are $P$-values for $(p_1,p_2)$ companies
from different two sections of NYSE: public utility section and capital goods section,
each of which has $n$ daily stock returns during the period 2000.1.1--2002.1.1. The number
of repeated experiments is $100$. All the closed stock prices are from WRDS database.
No. of Exp. is the number of experiments whose $P$-values are in the corresponding interval.}
\end{table}

%t11 #&#
\begin{table}%[h]
\tablewidth=242pt
\caption{$P$-values for $(p_1,p_2)$ companies from finance section and healthcare section of NYSE}\label{tbcc}
\begin{tabular*}{\tablewidth}{@{\extracolsep{\fill}}@{}lcc@{}}
\hline
& \multicolumn{2}{c@{}}{\textbf{No. of exp.}}\\[-6pt]
& \multicolumn{2}{c@{}}{\hrulefill}\\
\multicolumn{1}{@{}r}{$\bolds{P}$\textbf{-values:}} &  $\bolds{(p_1, p_2, n)}$ &  $\bolds{(p_1, p_2, n)}$\\
$\bolds{P}$\textbf{-value interval}  &   $\bolds{(10, 15, 20)}$ &  $\bolds{(15, 20, 25)}$\\
\hline
  $[0, 0.05]$   & 90 & 92 \\
  $[0.05, 0.1]$ & 4  & 5 \\
  $[0.1, 0.2]$  & 5 &  1 \\
  $[0.2, 0.3]$  & 1  & 2  \\
  $[0.3, 0.4]$  & 0  & 0 \\
  $[0.4, 0.5]$  & 0 &  0\\
  $[0.6, 0.7]$  & 0 &  0 \\
  $[0.8, 0.9]$  & 0 &  0 \\
  $[0.9, 1]$    & 0 &  0 \\
  \hline
\end{tabular*}
\tabnotetext[]{ta11}{These are $P$-values for $(p_1,p_2)$ companies from different two sections of NYSE: finance section and healthcare section,  each of which has $n$ daily stock returns during the period 2000.1.1--2002.1.1. The number of repeated experiments is $100$. All the closed stock prices are from WRDS database. No. of Exp. is the number of experiments whose $P$-values are in the corresponding interval.}
\end{table}

%s8 #&#
\section{Acknowledgement}\label{sec8}
The authors would like to thank the Editor, an Associate Editor and the
referees for their constructive comments and suggestions which
significantly improved this paper.

\begin{supplement}[id=suppA]
\stitle{Supplement to ``Independence test for high dimensional data based on 
regularized canonical correlation coefficients''\\}
\slink[doi]{10.1214/14-AOS1284SUPP} %[doi,text={...}] - jei reikia
%suskaldyti doi
\sdatatype{.pdf}
\sfilename{aos1284\_supp.pdf}
\sdescription{The supplementary material is divided into Appendices A~and~B.
Some useful lemmas, and proofs of all theorems and Proposition \ref{pr4}--\ref{pr5} are given in
Appendix~A while one theorem related to CLT of a sample covariance matrix plus a perturbation matrix is provided in Appendix~B.}
\end{supplement}

%===========bibliography===========================

% imsref loaded by linak, 2014-12-17 11:10:40
%
% imsref loaded by linak, 2014-12-23 14:06:47
% imsref loaded by linak, 2014-12-23 14:09:21

\printaddresses
\end{document}